\documentclass[10pt]{article}
\usepackage{amssymb,amsmath,latexsym}
\allowdisplaybreaks

\begin{document}

\begin{doublespace}

\def\1{{\bf 1}}
\def\ind{{\bf 1}}
\def\nn{\nonumber}

\def\sA {{\cal A}} \def\sB {{\cal B}} \def\sC {{\cal C}}
\def\sD {{\cal D}} \def\sE {{\cal E}} \def\sF {{\cal F}}
\def\sG {{\cal G}} \def\sH {{\cal H}} \def\sI {{\cal I}}
\def\sJ {{\cal J}} \def\sK {{\cal K}} \def\sL {{\cal L}}
\def\sM {{\cal M}} \def\sN {{\cal N}} \def\sO {{\cal O}}
\def\sP {{\cal P}} \def\sQ {{\cal Q}} \def\sR {{\cal R}}
\def\sS {{\cal S}} \def\sT {{\cal T}} \def\sU {{\cal U}}
\def\sV {{\cal V}} \def\sW {{\cal W}} \def\sX {{\cal X}}
\def\sY {{\cal Y}} \def\sZ {{\cal Z}}

\def\bA {{\mathbb A}} \def\bB {{\mathbb B}} \def\bC {{\mathbb C}}
\def\bD {{\mathbb D}} \def\bE {{\mathbb E}} \def\bF {{\mathbb F}}
\def\bG {{\mathbb G}} \def\bH {{\mathbb H}} \def\bI {{\mathbb I}}
\def\bJ {{\mathbb J}} \def\bK {{\mathbb K}} \def\bL {{\mathbb L}}
\def\bM {{\mathbb M}} \def\bN {{\mathbb N}} \def\bO {{\mathbb O}}
\def\bP {{\mathbb P}} \def\bQ {{\mathbb Q}} \def\bR {{\mathbb R}}
\def\bS {{\mathbb S}} \def\bT {{\mathbb T}} \def\bU {{\mathbb U}}
\def\bV {{\mathbb V}} \def\bW {{\mathbb W}} \def\bX {{\mathbb X}}
\def\bY {{\mathbb Y}} \def\bZ {{\mathbb Z}}
\def\R {{\mathbb R}} \def\RR {{\mathbb R}} \def\H {{\mathbb H}}
\def\n{{\bf n}} \def\Z {{\mathbb Z}}

\newcommand{\expr}[1]{\left( #1 \right)}
\newcommand{\cl}[1]{\overline{#1}}
\newtheorem{thm}{Theorem}[section]
\newtheorem{lemma}[thm]{Lemma}
\newtheorem{defn}[thm]{Definition}
\newtheorem{prop}[thm]{Proposition}
\newtheorem{corollary}[thm]{Corollary}
\newtheorem{remark}[thm]{Remark}
\newtheorem{example}[thm]{Example}
\numberwithin{equation}{section}
\def\ee{\varepsilon}
\def\qed{{\hfill $\Box$ \bigskip}}
\def\NN{{\mathcal N}}
\def\AA{{\mathcal A}}
\def\MM{{\mathcal M}}
\def\BB{{\mathcal B}}
\def\CC{{\mathcal C}}
\def\LL{{\mathcal L}}
\def\DD{{\mathcal D}}
\def\FF{{\mathcal F}}
\def\EE{{\mathcal E}}
\def\QQ{{\mathcal Q}}
\def\SS{{\mathcal S}}
\def\RR{{\mathbb R}}
\def\R{{\mathbb R}}
\def\L{{\bf L}}
\def\K{{\bf K}}
\def\S{{\bf S}}
\def\A{{\bf A}}
\def\E{{\mathbb E}}
\def\F{{\bf F}}
\def\P{{\mathbb P}}
\def\N{{\mathbb N}}
\def\eps{\varepsilon}
\def\wh{\widehat}
\def\wt{\widetilde}
\def\pf{\noindent{\bf Proof.} }
\def\pff{\noindent{\bf Proof} }
\def\cp{\mathrm{Cap}}

\title{\Large \bf Minimal thinness with respect to symmetric L\'evy processes}

\author{{\bf Panki Kim}\thanks{This work was supported by the National Research Foundation of
Korea(NRF) grant funded by the Korea government(MEST) (NRF-2013R1A2A2A01004822)
}
\quad {\bf Renming Song\thanks{Research supported in part by a grant from
the Simons Foundation (208236)}} \quad and
\quad {\bf Zoran Vondra\v{c}ek}
\thanks{Research supported in part by the Croatian Science Foundation under the project 3526}
}

\date{}

\maketitle

\begin{abstract}
Minimal thinness is a notion that describes the smallness of a set
at a boundary point. In this paper, we provide  tests for
minimal thinness at finite and infinite minimal Martin boundary points
for a large class of purely discontinuous
symmetric L\'evy processes.
\end{abstract}



%%%%%%%%%%%%%%%%%%%%%%%%%%%%%%%%%%%%%%%%%%%%%%%%%%%%%%%%%%%%%%%%%%%%%%%%%%%%%%%%%%%%%%%%%%%%%%%%%%%%%%%%%%%%%%%%%%%%%%%%%%%%%%%%%%
%%%%%%%%%%%%%%%%%%%%%%%%%%%%                          Introduction               %%%%%%%%%%%%%%%%%%%%%%%%%%%%%%%%%%%%%%%%%%%%%%%%%
%%%%%%%%%%%%%%%%%%%%%%%%%%%%%%%%%%%%%%%%%%%%%%%%%%%%%%%%%%%%%%%%%%%%%%%%%%%%%%%%%%%%%%%%%%%%%%%%%%%%%%%%%%%%%%%%%%%%%%%%%%%%%%%%%%

\section{Introduction}\label{int}
Minimal thinness
is a notion that describes the smallness of a set at a boundary point.
Minimal thinness in the half-space was introduced by Lelong-Ferrand in \cite{LF},
while minimal thinness in general open sets was developed by Na\"{\i}m in \cite{Nai};
for a more recent exposition see \cite[Chapter 9]{AG}. A probabilistic interpretation
in terms of Brownian motion was given by Doob, see e.g. \cite{Doo}.

A Wiener-type criterion for minimal thinness of a subset of the half-space
(using a Green energy) has already appeared in \cite{LF}. A refined version of
such a criterion (using the ordinary capacity) was proved in \cite{Ess}.
A general version of the Wiener-type criterion for minimal thinness in NTA
domains was established by Aikawa in \cite{A} using a powerful concept of
quasi-additivity of capacity. In case of a smooth domain, Aikawa's version of
the criterion implies several results obtained earlier, cf.~\cite{Beu, Dah, Maz, Sj77}.
A good exposition of this theory can be found in \cite[Part II, 7]{AE}.

All of these results have been proved in the context of classical potential theory
related to the Laplacian, or probabilistically, to Brownian motion.
Even though the concept of minimal thinness for Hunt processes admitting a dual process (and satisfying
an additional hypothesis)
was studied by F\"ollmer \cite{Fol},
concrete criteria for minimal thinness with respect to certain integro-differential
operators (i.e., certain L\'evy processes)
in the half space have been obtained only recently in \cite{KSV6}.
To be more precise,
in \cite{KSV6} the underlying process $X$ belongs to a class of subordinate Brownian motions,
where the Laplace exponents of the corresponding subordinators are complete
Bernstein functions satisfying a certain condition at infinity. The first
result of \cite{KSV6} was a necessary condition for minimal thinness of
a Borel subset $E$ of the half-space $\H\subset \R^d$, $d\ge 2$:
 If $E$ is minimally thin in $\H$ (with respect to the process $X$)
at the point $z=0$, then
$$
\int_{E\cap B(0,1)} |x|^{-d}\, dx <\infty\, .
$$
Here, and in the sequel, $B(z,r)$ denotes the open ball centered at $z\in \R^d$
with radius $r>0$. In the classical case this was proved in \cite{Beu}
for $d=2$ and in \cite{Dah} for $d\ge 3$. The method applied in
\cite{KSV6} was based on a result of  Sj\"ogren, \cite[Theorem 2]{Sj2}.
The second result of \cite{KSV6} was a criterion for minimal thinness in
$\H$ of a set under the graph of a Lipschitz function which in the classical
case is due to Burdzy, see \cite{Bur} and \cite{Gar}.

The goal of this paper is to generalize the results from \cite{KSV6} in several directions.
We will always assume that $d\ge 2$.
We  work with a broader class of purely discontinuous
symmetric transient L\'evy processes and
prove a version of Aikawa's Wiener-type criterion for
minimal thinness of a subset of a (not necessarily bounded) $\kappa$-fat open
set at any finite (minimal Martin) boundary point. By specializing to
$C^{1,1}$ open sets, we  get an integral criterion for minimal thinness
in the spirit of \cite{Beu, Dah}. Moreover,
in case the processes  satisfy an additional
assumption governing the global behavior, we obtain criteria for minimal thinness
of a subset of half-space-like open sets at infinity. In the classical case of the
Laplacian, such results are direct consequences of the corresponding finite boundary point results
by use of the inversion with respect to a sphere and the Kelvin transform.
In the case we study, this is much more delicate, since the method of Kelvin
transform is not available.

Let us describe the results of the paper in more detail. We start with a
description of the setup of this paper.

We assume that $r\mapsto j(r)$ is a strictly positive
non-increasing function
on $(0, \infty)$ satisfying
\begin{equation}\label{e:j-decay}
j(r)\le cj(r+1) \qquad \hbox{for } r\ge 1,
\end{equation}
and  $X$ is a purely discontinuous symmetric transient L\'evy process in $\R^d$
with L\'evy exponent $\Psi_X(\xi)$ so that
$$
\E_x\left[e^{i\xi\cdot(X_t-X_0)}\right]=e^{-t\Psi_X(\xi)}, \quad \quad t>0, x\in \R^d, \xi\in \R^d.
$$
We assume that the L\'evy measure of $X$ has a
density $J_X$ such that
\begin{equation}\label{e:psi1}
\gamma^{-1} j(|y|)\le J_X(y) \le \gamma j(|y|), \quad \mbox{for all } y\in \R^d\, ,
\end{equation}
for some $\gamma\ge 1$.
Since
$
\int_0^\infty j(r) (1\wedge r^2) r^{d-1}dr < \infty$ by \eqref{e:psi1},
the function $x \to j(|x|)$ is the L\'evy density of  an isotropic unimodal L\'evy process whose characteristic exponent is
$
\Psi(|\xi|)= \int_{\R^d}(1-\cos(\xi\cdot y))j(|y|)dy.
$
The L\'evy exponent $\Psi_X$ can be written as
$
\Psi_X(\xi)= \int_{\R^d}(1-\cos(\xi\cdot y))J_X(y)dy
$
and, clearly by \eqref{e:psi1}, it satisfies
\begin{equation}\label{e:psi21}
\gamma^{-1} \Psi(|\xi|)\le \Psi_X(\xi) \le \gamma \Psi(|\xi|), \quad \mbox{for all } \xi\in \R^d\, .
\end{equation}
The function $\Psi$ may be not increasing.
However, if we put $\Psi^*(r):= \sup_{s \le r} \Psi(s)$, then, by
\cite[Proposition 2]{BGR1}, we have
$$
\Psi(r) \le\Psi^*(r) \le \pi^2 \Psi(r).
$$
Thus by \eqref{e:psi21},
\begin{equation}\label{e:psi3}
(\pi^2\gamma)^{-1}  \Psi^*(|\xi|)\le \Psi_X(\xi) \le \gamma \Psi^*(|\xi|), \quad \mbox{for all } \xi  \in \R^d\, .
\end{equation}
Moreover,
\begin{equation}\label{e:lemma2.1}
{\Psi^*(\lambda t)}{} \le 2(1+\lambda^2) \Psi^*(t)\quad \text{ for every } t >0 \text{ and }\lambda \geq 1,
\end{equation}
(see \cite[Lemma 1]{G}).

We will always assume that $\Psi$ satisfies the following scaling condition at infinity:

\medskip
\noindent
({\bf H1}):
There exist constants $0<\delta_1\le \delta_2 <1$ and $a_1, a_2>0$  such that
\begin{equation}\label{e:H1n}
a_1\lambda^{2\delta_1} \Psi(t) \le \Psi(\lambda t) \le a_2 \lambda^{2\delta_2}
\Psi(t), \quad \lambda \ge 1, t \ge 1\, .
\end{equation}
Then by  \cite[(15) and Corollary 22]{BGR1},
for every $R>0$, there exists $c=c(R)>1$ such that
\begin{equation}\label{e:asmpbofjat0n}
c^{-1}\frac{\Psi(r^{-1})}{r^d} \le j(r)\le c \frac{\Psi(r^{-1})}{r^d}
\quad \hbox{for } r\in (0, R].
\end{equation}
Note that the class of L\'evy processes described above contains
the purely discontinuous unimodal L\'evy processes
dealt with in \cite{BGR1, BGR3}.
The condition ({\bf H1})
governs the small time and small space
behavior of the process $X$. Thus it is sometimes referred to as a {\em local} condition.

\medskip
{\it
In this paper,  we will always assume that the  condition {\bf (H1)} is satisfied and
$X$ is a purely discontinuous symmetric transient L\'evy process with L\'evy density $J_X$
satisfying \eqref{e:psi1}.}
\medskip

To study minimal thinness at infinity, we need to add another scaling condition on
$\Psi$  near the origin:

\medskip
\noindent
({\bf H2}):
There exist constants $0<\delta_3\le \delta_4 <1$ and $a_3, a_4>0$  such that
\begin{equation}\label{e:H2}
a_3\lambda^{2\delta_4} \Psi(t) \le \Psi(\lambda t) \le a_4 \lambda^{2\delta_3}
\Psi(t), \quad \lambda \le 1, t \le 1\, .
\end{equation}

Since $d\ge 2$, ({\bf H2}) implies that $X$ is transient.
Condition ({\bf H2}) governs the large time and large space
behavior of $X$ and so it is sometimes referred to as a {\em global} condition.

We will impose
the condition ({\bf H2}) only when we discuss
minimal thinness at infinity and we will explicitly mention this assumption when stating the results or at the beginning of the section.

Let $({\cal E}, {\cal F})$ be the Dirichlet form of $X$ on $L^2(\R^d, dx)$.
It is known that $({\cal E}, {\cal F})$
is a regular Dirichlet form on $L^2(\R^d, dx)$
and ${\cal F}=\overline{C^\infty_c(\R^d)}^{\, {\cal E}_1}
= \{f\in L^2(\R^d, dx): {\cal E}(f, f)<\infty\}$,
where ${\cal E}_1(u, u)={\cal E}(u, u) + \int_{\R^d}u^2(x)dx$.
 Moreover, for $u\in {\cal F}$,
$$
{\cal E}(u, u)=\int_{\R^d\times \R^d}(u(x)-u(y))^2
J_X(x-y)dxdy.
$$

For any open $D\subset \R^d$, we use $X^D=(X^D_t, \P_x)$ to denote the subprocess of $X$ killed upon
exiting $D$. The Dirichlet form of $X^D$ is given by $({\cal E}, {\cal F}_D)$, where
$$
{\cal F}_D=\{u\in {\cal F}: u=0  \mbox{ on $D^c$ except for
a set of zero capacity}\}.
$$

The Hardy inequality is one of the main ingredients in Aikawa's construction of a
measure comparable to the capacity, which is fundamental in proving quasi-additivity
of capacity
(see \cite{A, AE}).
 We introduce a local Hardy inequality for the Dirichlet form $(\EE, \FF_D)$
in the next definition and show in Section \ref{qac} that it holds under natural conditions on the open set $D$.

\begin{defn}\label{def:LHI}
We say that $({\cal E}, {\cal F}_D)$ satisfies
the local Hardy inequality at $z \in \partial D$ (with a localization constant $r_0$)
if there exist $c>0$ and  $r_0>0$ such that
$$
\EE (v,v)\ge c \int_{B(z, r_0) \cap D} v^2(x)\Psi(\delta_D(x)^{-1})\, dx\, , \quad v\in \FF_D\, .
$$
\end{defn}

We recall now the definition of $\kappa$-fat open set and introduce the necessary notation.
\begin{defn}\label{def:UB}
Let $0<\kappa\leq 1$. We say that an open set $D$ is $\kappa$-fat at $z \in \partial D$ if there is a localization radius $R>0$ such that for all $r\in (0,R]$ there exists a (non-tangential) point $A_r(z)\in D\cap B(z,r)$ such that the ball $B(A_r(z) , \kappa r) \subset D\cap B(z,r)$. We say that an open set is $\kappa$-fat with localization radius $R$ if $D$ is $\kappa$-fat
at all $z \in \partial D$ with localization radius $R$.
\end{defn}

Without loss of generality, we will assume that $R\le 1/2$ and $\kappa\le 1/4$.

The first main result of this paper is the following Aikawa's version of the Wiener-type
criterion for minimal thinness.
For any open set $D\subset \R^d$, we use $G_D$ to denote the Green function of $X^D$.
See Definition \ref{def:cMthin} for
the definition of minimal thinness in $D$ with respect to $X$.
\begin{thm}\label{t:dahlberg}
Assume that $D\subset \R^d$ is $\kappa$-fat with  localization radius $R$
and that $({\cal E}, {\cal F}_D)$ satisfies
the local Hardy inequality with a localization constant $r_0$ at $z \in \partial D$.
Fix a point $x_0 \in D$ with $\kappa R < \delta_D(x_0) <R$.

\noindent
(1)
If a Borel set $E\subset D$ is  minimally thin in $D$  at $z$ with respect to $X$, then
$$
\int_{E\cap B(z,(\kappa R/4) \wedge r_0)} \left( \frac{G_D(x,x_0)}
{G_D(A_{|x-z|} (z),x_0 )}\right)^2 \frac{\Psi(\delta_D(x)^{-1})}
{\Psi(|x-z|^{-1})}\,\frac{dx}{|x-z|^d}< \infty\, .
$$

\noindent
(2)
Conversely, if $E$ is the union of a subfamily of Whitney cubes of $D$ and is
not minimally thin in $D$  at $z$ with respect to $X$, then
$$
\int_{E\cap B(z,(\kappa R/4) \wedge r_0)} \left( \frac{G_D(x,x_0)}
{G_D(A_{|x-z|} (z),x_0)}\right)^2 \frac{\Psi(\delta_D(x)^{-1})}{\Psi(|x-z|^{-1})}\,
\frac{dx}{|x-z|^d}=\infty\, .
$$
\end{thm}

When $D$ is a half space, or when
$D$ is a $C^{1,1}$ open set
and $X$ is a purely discontinuous unimodal L\'evy processes,
we have an explicit form of the integral test.
We first recall the definition of a $C^{1,1}$ open set.
\begin{defn}
An open set $D$ in $\bR^d$ is said to be a (uniform)
$C^{1,1}$ open set if there exist a localization radius $R>0$ and
a constant $\Lambda>0$ such that for every $z\in\partial D$, there
exist a $C^{1,1}$-function $\psi=\psi_z: \bR^{d-1}\to \bR$
satisfying $\psi (0)= 0$,
$\nabla \psi (0)=(0, \dots, 0)$, $\| \nabla \psi  \|_\infty \leq \Lambda$,
$| \nabla \psi (x)-\nabla \psi (w)| \leq \Lambda |x-w|$,
and an orthonormal coordinate system $CS_z$ with its origin at $z$ such that
$$
B(z, R_1)\cap D=\{ y= (\wt y, \, y_d) \mbox{ in } CS_z: |y|< R,
y_d > \psi (\wt y) \}.
$$
The pair $(R, \Lambda)$ is called the characteristics of the $C^{1,1}$ open set $D$.
\end{defn}

A $C^{1,1}$ open set $D$ with characteristics $(R, \Lambda)$ can be unbounded and disconnected;
the distance between two distinct components of $D$ is at least $R$.

Recall that an open set $D$ is said to satisfy the interior and exterior balls conditions with radius
$R_1$ if for every $z\in \partial D$, there exist $x\in D$ and $y\in \overline{D}^c$ such that
${\rm dist}(x, \partial D)=R_1$, ${\rm dist}(y, \partial D)=R_1$,
$B(x, R_1)\subset D$ and $B(y, R_1)\subset  \overline{D}^c$.

It is known, see \cite[Definition 2.1 and Lemma 2.2]{AKSZ}, that an open set
$D$ is a $C^{1,1}$ open set if and only if it satisfies the interior and exterior ball conditions.
By taking $R$ smaller if necessary, we will always assume a  $C^{1,1}$ open set
with characteristics $(R, \Lambda)$ satisfies the interior and exterior balls conditions with radius
$R$.

\begin{corollary}\label{c:dahlberg}
Suppose that either (i) $D$ is a half space; or (ii) $D\subset \R^d$ is a $C^{1,1}$ open set and $\gamma=1$
in \eqref{e:psi1}.
Assume that  $E$ is a Borel subset of $D$.

\noindent
(1)
If $E$ is  minimally thin in $D$  at $z\in\partial D$ with respect to $X$, then
$$
\int_{E\cap B(z,1)} |x-z|^{-d}\, dx < \infty\, .
$$

\noindent
(2)
Conversely, if $E$ is the union of a subfamily of Whitney cubes of $D$ and is
not minimally thin in $D$  at $z\in\partial D$ with respect to $X$, then
$$
\int_{E\cap B(z,1)} |x-z|^{-d}\, dx =\infty\, .
$$
\end{corollary}

We sometimes write a point $z=(z_1, \dots, z_d)\in \bR^d$ as
$(\wt z, z_d)$ with $\wt z \in \bR^{d-1}$.
Throughout this paper, $\bH_b$ stands for the set $\{x=(\wt x, x_d)\in \R^d: x_d>b\}$. We will denote the upper half space $\bH_0$ by $\bH$.

An open set $D$ is said to be
half-space-like if, after isometry, there exist two real numbers
$b_1\le b_2$ such that $\bH_{b_2}\subset D\subset \bH_{b_1}$.
Without loss of generality, whenever we deal with a half-space-like open set $D$,
we will always assume that
$
\bH_{1}\subset D\subset \bH.
$

Now we state our results on minimal thinness at infinity.
In Section \ref{mti} we will first extend the main result of \cite{KSV9} to
purely discontinuous unimodal L\'evy processes so that,
for a large class of unbounded open sets including half-space-like open sets,
the infinite part of the (minimal) Martin boundary consists of a single point.
We call such a point infinity and denote it by $\infty$.

\begin{defn}\label{def:HI}
 We say $({\cal E}, {\cal F}_D)$ satisfies the Hardy inequality if there exists $c>0$ such that
$$
\EE(v,v)\ge c \int_D v^2(x)\Psi(\delta_D(x)^{-1})\, dx\, , \quad v\in \FF_D\, .
$$
\end{defn}

Here is the second main result of the paper.

\begin{thm}\label{t:dahlberg-infinity}
Suppose that  {\rm ({\bf H2})} holds and  $\gamma=1 $ in \eqref{e:psi1}.
Assume
that $D\subset \R^d$ is a half-space-like open set  and that
$({\cal E}, {\cal F}_D)$ satisfies the Hardy inequality.
 Let $E$ be a Borel subset of $D$ and
$x_0=(\wt 0, 5)$.

 \noindent
(1)  If $E$ is  minimally thin in $D$  at infinity with respect to $X$, then
$$
\int_{E\cap B(0,1)^c} |x|^{d} G_D(x,x_0)^2 \Psi(\delta_D(x)^{-1}) \, dx < \infty\, .
$$

 \noindent
(2)
Conversely, if $E$ is the union of a subfamily of Whitney cubes of $D$ and is
not minimally thin in $D$  at  infinity with respect to $X$, then
$$
\int_{E\cap B(0,1)^c} |x|^{d} G_D(x,x_0)^2 \Psi(\delta_D(x)^{-1})\, dx =\infty\, .
$$
\end{thm}

Again, when $D$ is a half-space-like $C^{1,1}$ open set, we get the following corollary.
\begin{corollary}\label{c:dahlberg-infinity}
Suppose that  {\rm ({\bf H2})} holds and $\gamma=1$ in \eqref{e:psi1}.
Assume
that $D\subset \R^d$ is a half-space-like $C^{1,1}$ open set and that $E$ is a Borel subset of $D$.

 \noindent
(1)  If $E$ is  minimally thin in $D$  at infinity with respect to $X$, then
$$
\int_{E\cap B(0,1)^c} |x|^{-d}\, dx < \infty\, .
$$

 \noindent
(2)
Conversely, if $E$ is the union of a subfamily of Whitney cubes of $D$ and is
not minimally thin in $D$  at  infinity with respect to $X$, then
$$
\int_{E\cap B(0,1)^c} |x|^{-d}\, dx =\infty\, .
$$
\end{corollary}

In order to prove these results we need various potential-theoretic results for the
underlying process such as Harnack inequality, boundary Harnack principle,
sharp estimates of the Green function and the Martin kernel in $D\subset \R^d$,
and identification of the Martin boundary of $D$ with the Euclidean boundary.
All of these results have been established previously, some of them quite recently,
under various conditions on the process $X$ and the open set $D$.
The main novelty is
that local results for possibly unbounded open sets are obtained only under
local conditions on the underlying precess $X$ -- a fact that leads to significant
technical difficulties.
Therefore we start the paper with three preliminary sections
that establish all necessary results.
In Section \ref{gfe}
we first recall some previous results from \cite{KSV10}.
The main new result is Theorem \ref{t:green} where we prove sharp local estimates of
the Green function of $X^D$ in case $D$ is a (not necessarily bounded) $\kappa$-fat
open set.
It is proved in  \cite{KSV10} that the finite part of the (minimal) Martin boundary of any $\kappa$-fat set $D$
coincides with the Euclidean boundary of $D$, see \cite[Theorem 3.13]{KSV10}.
The main result of Section \ref{mke} is Theorem \ref{t:mke} which gives sharp estimates of the
 Martin kernel at a finite boundary point of a $\kappa$-fat open set.
In Section \ref{H1-H2},
we assume {\rm ({\bf H2})} holds and
extend
some results previously known for subordinate Brownian motions.

In Section \ref{qac}, we
will discuss both local results and  global results (under the condition ({\bf H2})).
In  that section, we study quasi-additivity of capacity with respect to a
Whitney decomposition of $D$. Here we closely follow the method from \cite{A},
but use more delicate estimates for the underlying L\'evy process.
The main novelty here is that we prove local quasi-additivity only under
local assumptions on the process $X$, see Proposition \ref{p:quasi-additivity}(1).
One of the main ingredients in proving quasi-additivity is a construction of a
measure comparable to capacity. Here one needs a Hardy-type inequality for
the associated Dirichlet form. We assume the (local) Hardy inequality and
at the end of the section give some sufficient conditions for this inequality.

Section \ref{mtf}
is devoted to the proof of Theorem \ref{t:dahlberg}.
After recalling the definition of minimal thinness and giving a general criterion,
we establish in Lemma \ref{l:key-step-1} the main technical tool for proving
the Wiener-type criterion for minimal thinness at a finite boundary point given
in Proposition \ref{p:minthin-criterion-1}.  Arguments leading from this criterion to
its Aikawa's version given in Proposition \ref{p:aikawa-thinness} are analogous to
those of \cite[Part II, 7]{AE} and rely on the (local) quasi-additivity.  The proof of
the main Theorem \ref{t:dahlberg} is a consequence of Aikawa's criterion and
the existence of a comparable measure.
In the case when $X$ is a unimodal L\'evy process  and $D$ a $C^{1,1}$-open set,
explicit boundary behaviors of the mean exit times in terms of the distance
to the boundary lead to
Corollary \ref{c:dahlberg}, see also the proof of
Corollary  \ref{c:aikawa-thinness} --  Aikawa's Wiener-criterion for $C^{1,1}$ open set.

In Section \ref{mti}
we assume  that {\rm ({\bf H2})} holds and that $X$ is a unimodal L\'evy process .  In this section we
study minimal thinness at infinity under global assumptions
on the underlying process and prove Theorem \ref{t:dahlberg-infinity}. The proofs,
although similar to these from the previous section, contain non-trivial modifications
(in particular the main technical Lemma \ref{l:Gfeinfinty}) and are given in full.
The starting point of the section is Theorem \ref{t:mbatinfty} where we
extend a
recent result from \cite{KSV9} stating that  the (minimal) Martin boundary of an open set
which is $\kappa$-fat at infinity consists of precisely one point. Besides half-space-like
open sets, infinite cones are another example of unbounded sets which are
$\kappa$-fat at infinity. Minimal thinness at infinity for infinite cones seems to
be more delicate, even in the classical case, see \cite[Theorem 1]{MY}.
That is why we restrict our
consideration to half-space-like open sets.

Finally, in Section \ref{graph} we  study
the question of minimal thinness of a set
below the graph of a Lipschitz function, both at a finite and infinite boundary point.
In case of minimal thinness at a finite boundary point we state in Proposition
\ref{c:dahlberg2} a Burdzy's type criterion which generalizes \cite[Theorem 4.4]{KSV6}.
As an application of Theorem \ref{t:dahlberg-infinity} and Corollary \ref{c:dahlberg-infinity}
we prove the main result of
 Section \ref{graph} --  a criterion for minimal thinness at infinity under the graph of
a Lipschitz function, see Theorem \ref{t:main2}.

We conclude this introduction by setting up some notation and conventions.
We use ``$:=$" to denote a definition, which
is  read as ``is defined to be"; we denote $a \wedge b := \min \{ a, b\}$,
$a \vee b := \max \{ a, b\}$;
we often denote point $z=(z_1, \dots, z_d)\in \bR^d$ as
$(\wt z, z_d)$ with $\wt z \in \bR^{d-1}$;
 we
denote by $B(x, r)$ the open
ball centered at $x\in \bR^d$ with radius $r>0$;
for any two positive functions $f$ and $g$,
$f\asymp g$ means that there is a positive constant $c\geq 1$
so that $c^{-1}\, g \leq f \leq c\, g$ on their common domain of
definition;
for any Borel subset $E\subset\bR^d$ and $x\in E$,
$\mathrm{diam}(E)$ stands for the diameter of $E$,
$|E|$ stands for the Lebesgue measure of $E$ in $\R^d$,
int$(E)$ stands for the interior of $E$
and
$\delta_E(x)$ stands for the Euclidean distance between
$x$ and $E^c$; $\N$ is the set of nonnegative integers.
The values of
the constants  $R, \delta_1, \delta_2, \delta_3, \delta_4 ,C_1, C_2, \dots$ remain the same
throughout this paper, while $c, c_0, c_1, c_2, \ldots$  represent constants
whose values are unimportant and may change. All constants are positive finite numbers.
The labeling of the constants $c_0, c_1, c_2, \ldots$ starts anew in the statement
and proof of each result.
The dependence of constant $c$ on dimension $d$ is not mentioned explicitly.

%%%%%%%%%%%%%%%%%%%%%%%%%%%%%%%%%%%%%%%%%%%%%%%%%%%%%%%%%%%%%%%%%%%%%%%%%%%%%%%%%%%%%%%%%%%%%%%%%%%%%%%%%%%%%%%%%%%%%%%%%%%%%%%%%%
%%%%%%%%%%%%%%%%%%%%%%%%%%%%                  Green function estimates           %%%%%%%%%%%%%%%%%%%%%%%%%%%%%%%%%%%%%%%%%%%%%%%%%
%%%%%%%%%%%%%%%%%%%%%%%%%%%%%%%%%%%%%%%%%%%%%%%%%%%%%%%%%%%%%%%%%%%%%%%%%%%%%%%%%%%%%%%%%%%%%%%%%%%%%%%%%%%%%%%%%%%%%%%%%%%%%%%%%%

\section{Green function estimates}\label{gfe}

Throughout this paper, we always assume  that $j$ is a strictly positive
non-increasing function on $(0, \infty)$ satisfying \eqref{e:j-decay} such
that ${\bf (H1)}$ holds, and that
$X$ is a purely discontinuous symmetric transient  L\'evy
process with L\'evy exponent $\Psi_X(\xi)$ and  a L\'evy density $J_X$ satisfying \eqref{e:psi1}.

As a consequence of {\bf (H1)}, \eqref{e:psi3} and \cite[Proposition 28.1]{Sa} we
know that for any $t>0$, $X_t$ has a density $p_t(x, y)$ which is smooth.
We will use $G(x, y):=\int_0^\infty p_t(x,y) dt $ to denote the Green function of $X$.
Since $X$ is a L\'evy process, $G(x,y)$ depends on $x-y$ only. Moreover, by the symmetry assumption on $X$
 $G(x,y)=G(y,x)$. For simplicity we use $G(x)$ for $G(x,0)$.

Given  an open set $D\subset \R^d$, we define
$X^D_t(\omega)=X_t(\omega)$ if $t< \tau_D(\omega)$ and
$X^D_t(\omega)=\partial$ if $t\geq  \tau_D(\omega)$, where
$\partial$ is a cemetery state.
Throughout this paper, we use the convention that any function $f$ on $D$ is extended to the cemetery state by $f(\partial)=0$.
Since $J_X$ satisfies the assumption \cite[(1.6)]{CKK2},
by \cite[Theorem 3.1]{CKK2}, for every
open set $D$, $X^D_t$ has a H\"older  continuous density $p_D(t,x,y)$.
For any
open  set $D$ in $\R^d$, let $G_D(x,y)=\int_0^{\infty}p_D(t,x,y)\, dt$
be the Green function of $X^D$.
The function $G_D$ is jointly lower semi-continuous on $D\times D$, see Remark \ref{r:lsc}.

We first recall the definitions of harmonic functions
with respect to $X$ and $X^D$.
\begin{defn}\label{def:har1}
Let $D$ be an open subset of $\R^d$.
A nonnegative function $u$ on $\R^d$ is said to be
(1) harmonic in $D$ with respect to $X$ if
$u(x)= \E_x\left[u(X_{\tau_{B}})\right]$ for each $x\in B$
and  every open set $B$ whose closure is a compact
subset of $D$;
(2)
regular harmonic in $D$ with respect to $X$ if
for each $x \in D$,
$
u(x)= \E_x\left[u(X_{\tau_{D}}), \tau_D<\infty\right].
$
\end{defn}

\begin{defn}\label{def:harkilled}
Let $D$ be an open subset of $\R^d$.
A nonnegative function $u$ on $D$ is said to be
harmonic with respect to $X^D$ if
$
u(x)= \E_x\left[u(
X^D_{\tau_{
U}})\right]$ for every $x\in
U$ and every open set $
U$ whose closure is a compact
subset of $D$.
\end{defn}

Obviously, if $u$ is harmonic with respect to $X^D$, then the function equal to
$u$ in $D$ and zero outside $D$ is harmonic with respect to $X$ in $D$.
All nonnegative functions that are harmonic in
$D$ with respect to $X$ are continuous, see \cite{KSV10}.

For notational convenience,  we define
$$
\Phi(r)=\frac1{\Psi^*(r^{-1})}, \quad r>0.
$$
The right continuous inverse function of $\Phi$ will be denoted by the usual notation $\Phi^{-1}(r)$.

The following two results are proved in \cite{KSV10}.

\begin{thm}[{\cite[(1.4), (2.1) and Theorem 2.19]{KSV10}}]\label{t:greenh}
For every $M \ge 1$  there exists
$C_1(M)=C_1(M, \Psi, \gamma)>0$
such that for all  $x\in B(0, M)$,
$$
 C_1(M)^{-1} \frac{\Phi(|x|)}{|x|^{d}} \, \le\, G(x)
\,\le \, C_1(M)
\frac{\Phi(|x|)}{|x|^{d}}\, .
$$
\end{thm}

\begin{lemma}[{\cite[(1.4), (2.1) and Lemma 2.12]{KSV10}}]\label{dec}
For every bounded open set $D$,
the Green function
$G_D(x, y)$ is finite and continuous off the diagonal of $D\times D$ and
there exists $c =c(\mathrm{diam}(D)$, $\Psi, \gamma)\ge 1$
such that for all $x, y\in D$,
$$
G_D(x, y) \le c \frac{\Phi(|x-y|)}{|x-y|^{d}}.
$$
\end{lemma}

Before we state the interior lower bound on the Green function, we first recall a result
from analysis (see \cite[Theorem 1, p. 167]{Stein}): Any open set
$D\subset\R^d$ is the union of a family $\{Q_j\}_{j \in \N}$ of closed cubes, with sides all parallel to the axes,
satisfying the following properties:
(i)  $\mathrm{int}(Q_j) \cap\, \mathrm{int} (Q_k)=\emptyset$;
(ii) for any $j$,
${\rm diam}(Q_j)\le {\rm dist}(Q_j, \partial D)\le 4  {\rm diam}(Q_j)$,
where ${\rm dist}(Q_j, \partial D)$ denotes the Euclidean distance between $Q_j$ and $\partial D$.
The family $\{Q_j\}_{j \in \N}$ above is called a Whitney decomposition of $D$ and the $Q_j$'s
are called Whitney
cubes (of $D$). We will use $x_j$ to denote the center of the cube $Q_j$.

\begin{lemma}\label{G:g3}
(1)
For every $L , T>0$, there exists $c=c(T, L, \Psi, \gamma)>0$ such
that for any bounded open set $D$ with $\mathrm{diam}(D) \le T$,
$x, y\in D$ with $|x-y| \le L ( \delta_D(x) \wedge \delta_D(y))$,
\begin{equation}\label{e:g3}
G_D(x,y) \ge c \frac{\Phi(|x-y|)}{|x-y|^d}.
\end{equation}

\noindent
(2) For every $M\ge 1$, every $L>0$ and any open set $D$,
there exists $c=c(M,L, \Psi,\gamma)>0$
such that for every Whitney decomposition $\{Q_j\}$ of $D$, every cube $Q_j$ such that
$\mathrm{diam}(Q_j)\le M$ and all $x,y\in Q_j$ with $|x-y|\le L(\delta_D(x)\wedge \delta_D(y))$,
$$
G_D(x,y) \ge c \frac{\Phi(|x-y|)}{|x-y|^d}.
$$
\end{lemma}

\pf
(1) This part is exactly \cite[Lemma 2.14]{KSV10}.

\noindent
(2) Fix $Q_j$ with $\mathrm{diam}(Q_j)\le M$. Recall that $x_j$ is the center of $Q_j$.
Let $\wt{D}:=D\cap B(x_j, 8M)$, so that $\mathrm{diam}(\wt{D})=8M$. Let $x,y\in Q_j$
be such that $|x-y|\le L(\delta_D(x)\wedge \delta_D(y))$. Since
$\mathrm{dist}(Q_j,\partial D)\le 4\mathrm{diam}(Q_j)\le 4M$,  we see that
$\delta_{\wt{D}}(x)=\delta_D(x)$ and $\delta_{\wt{D}}(y)=\delta_D(y)$.
Thus $|x-y|\le L(\delta_{\wt{D}}(x)\wedge \delta_{\wt{D}}(y))$. From part
(1)
we conclude that \eqref{e:g3} holds with
$c=c(M,L,\Psi,\gamma)$.
\qed

\begin{remark}\label{r:lsc}
By the domain monotonicity of Green functions and Lemmas \ref{dec} and \ref{G:g3}(i),
one can easily  see that
the function $G_D$ is jointly lower semi-continuous on $D\times D$.
In fact, $G_D$ is continuous in the extended sense at the diagonal:
$\lim_{(x,y)\to (x_0,x_0)}G_D(x,y)=G_D(x_0,x_0)=+\infty$ for any $x_0\in D$.
\end{remark}

We record a simple consequence of  {\bf (H1)}, which we will use several times:
There exists a positive constant $c_1>0$ such that
 for all positive $r,s,A$ with $As \le r \le 1$,
 \begin{align}\label{e:Gwd}
 \frac{\Phi(r)}{r^d} \le c_1  \left(  A^{-d} \vee  A^{-d+2\delta_1} \right)  \frac{\Phi(s)}{s^d}.
 \end{align}

In the remainder of this section, we assume that $D$ is a
$\kappa$-fat open set with localization radius $R$.
Without loss of generality we may assume that $R\le \frac{1}{10}$.
Recall that  for each $z \in \partial D$ and $r \in (0, R)$,
$A_r(z)$ is a point in
$ D \cap B(z,r)$ satisfying $ B(A_r(z),\kappa r)\subset D \cap B(z,r)$.
We also
recall that $G_D(\cdot, y)$ is
regular harmonic in $D\setminus \overline{B(y,\varepsilon)}$ for
every $\varepsilon >0$ and vanishes outside $D$.

\begin{lemma}\label{C:c_L} (Carleson's estimate)
There exists a constant $c=c(\Psi, \gamma, \kappa) >1$ such that for every
$z\in\partial D$,  $0<r  \le \kappa R/4$, $y \in D  \setminus \overline{B(z, 4r)}$,
\begin{equation}\label{e:CG_3}
G_{D}(x, y) \, \le \,c\,  G_{D} ( A_r(z), y), \quad x \in D \cap
B(z, r).
\end{equation}
\end{lemma}

\pf
By the boundary Harnack principle in \cite[Theorem 2.18(ii)]{KSV10},
it suffices to show that for $y \in D  \setminus \overline{B(z, 4r)}$,
 \begin{align}\label{e:c_L}
{G_{D}(x, A_{4r/\kappa}(z))}  \le c_1r^{-d} \Phi(r) \le c_2 {G_{D}(A_r(z), A_{4r/\kappa}(z))}, \quad x \in D \cap
B(z, r).
 \end{align}
Since $|x-A_{4r/\kappa}(z)| \ge
\delta_D(A_{4r/\kappa}(z))-\delta_D(x)\ge 4r-r=3r$,
the first inequality in \eqref{e:c_L} follows from  Theorem \ref{t:greenh} and \eqref{e:Gwd}.
 On the other hand,
since
$
|A_r(z)-A_{4r/\kappa}(z)| \le{8r}{\kappa}^{-1}\le8\kappa^{-2} (\delta_D(A_{4r/\kappa}(z)) \wedge
\delta_D(A_r(z))),$
the second inequality in \eqref{e:c_L} follows from
Lemma \ref{G:g3}(1) and \eqref{e:Gwd}.
The assertion of the lemma follows.
 \qed

Next we show a localization result for unbounded $\kappa$-fat set.
\begin{prop}\label{p:localization}
Let $D$ be an unbounded $\kappa$-fat set with localization radius $R>0$.
There exist $\kappa_1=\kappa_1(\kappa, R,d)\in (0,\kappa]$ and $R_1=R_1(\kappa,R,d)\in (0,R]$
such that for every $z_0\in \partial D$ there exists a $\kappa_1$-open set $D(z_0)$
 with localization radius $R_1$
satisfying $D\cap B(z_0, 1) \subset D(z_0) \subset D\cap B(z_0, 2)$.
\end{prop}

\pf
Recall that $R\le 1/4$ and $\kappa\le 1/2$.
We first note that by making $\kappa$ smaller if necessary,
we can assume that for any $r\le R$ and $z\in \overline{D}$, there exists $A_r(z)\in D$ such that $B(A_r(z), \kappa r)\subset D\cap B(z, r)$.
This means that the $\kappa$-fatness of $D$ at every boundary point $z\in \partial D$ implies that the $\kappa$-fatness property (with possibly smaller $\kappa$) holds  also true at every interior point $z\in D$. In particular, a non-tangential point $A_r(z)\in D\cap B(z,r)$ is well defined for every $z\in \overline{D}$.

Let $z_0\in \partial D$ and define
$$
D(z_0)=\left(D\cap B(z_0, 1)\right)\bigcup\left(\bigcup_{z\in \overline{D\cap B(z_0, 1)}}
\bigcup_{r\in (0, R]}B( A_r(z), \kappa r)\right).
$$
Clearly, $D\cap B(z_0, 1) \subset D(z_0) \subset D\cap B(z_0, 2)$.
We claim that $D(z_0)$ is $\kappa_1$-fat with localization radius $R$,
where $\kappa_1=\kappa/32$. Actually, we show that for every $w\in \overline{D(z_0)}$ and
every $s\in (0,R]$ there exists $\wt{A}_s(w)\in D(z_0)$ such that $B(\wt{A}_s(w), \kappa_1 s)\subset D(z_0)\cap B(w,s)$.

Let $w\in \overline{D\cap B(z_0,1)}$. By the very definition of $D(z_0)$, for each $s\in (0,R]$ we have that $B( A_s(w), \kappa s)\subset D(z_0)\cap B(w,s)$. Hence, we can take $\wt{A}_s(w)=A_s(w)$.

Suppose now that $w\in \overline{B(A_r(z),\kappa r)}$
 for some $z\in \overline{D\cap B(z_0,1)}$ and $r\in (0,R]$. We consider two cases: (i) $s\le 8r$, and (ii) $8r<s\le R$. In the case $s\le 8r$, first note that $2\kappa_1 s=(\kappa/16)s\le \kappa r/2$. Consider the line segment connecting $w$ and $A_r(z)$ and let $\wt{A}_s(w)$ be the point on this segment at the distance $2\kappa_1 s$ from $w$. If $v\in B(\wt{A}_s(w), \kappa_1 s)$, then
$$
|v-w|\le |v-\wt{A}_s(w)|+
|\wt{A}_s(w)-w|
\le \kappa_1 s+2\kappa_1 s\le \frac{3}{32}\kappa s<\frac{s}{2}\, ,
$$
and
$|v-A_r(z)|\le |w-A_r(z)|\le \kappa r\, .$
This proves that
$B(\wt{A}_s(w), \kappa_1 s)\subset B(A_r(z),\kappa r)\cap B(w,s/2)\subset D(z_0)\cap B(w,s/2)$.
If $8r<s\le R$, then $B(z, s/4)\subset B(w,s/2)$. Indeed, if $v\in B(z,s/4)$, then
$$
|v-w|\le |v-z|+|z-A_r(z)|+|A_r(z)-w|\le \frac{s}{4}+r+\kappa r \le \frac{s}{4}+\frac32 r<\frac{s}{2}\, .
$$
Take $\wt{A}_s(w)=A_{s/4}(z)$. Then $B(\wt{A}_s(w), \kappa_1 s)\subset B(A_{s/4}(z),\kappa s/4)\subset D(z_0)\cap B(z,s/4)\subset D(z_0)\cap B(w,s/2)$.

Finally, assume that $w$ is in the closure of $\bigcup_{z\in \overline{D\cap B(z_0, 1)}}
\bigcup_{r\in (0, R]}B( A_r(z), \kappa r)$. Then there exist sequences $(w_n)_{n\ge 1}$,
$(z_n)_{n\ge 1}$ in $\overline{D\cap B(z_0,1)}$ and $(r_n)_{n\ge 1}$ in ($0,R]$ such that
$w_n\in B(A_{r_n}(z_n), \kappa r_n)$ and $w=\lim_n w_n$. Let $s\le R$ and choose $n\ge 1$
so that $|w_n-w|\le s/4$.
By what has already been proved, there exists $\wt{A}_s(w_n)\in D(z_0)$
so that $B(\wt{A}_s(w_n), \kappa_1 s)\subset D(z_0)\cap B(w_n, s/2)$.
Define $\wt{A}_s(w)=\wt{A}_s(w_n)$. If $v\in B(\wt{A}_s(w), \kappa_1 s)$, then
$$
|v-w|\le |v-\wt{A}_s( w_n)|+|\wt{A}_s(w_n)-w_n|+|w_n-w|\le \kappa_1 s+\frac{s}{2}+\frac{s}{4}<s\, .
$$
Thus $B(\wt{A}_s(w), \kappa_1 s)\subset D(z_0)\cap B(w,s)$.
\qed

Without loss of generality (by choosing  $R$ and $\kappa$ smaller if necessary), we assume in the sequel that $R=R_1$ and $\kappa=\kappa_1$.

Fix $x_0 \in D$ with $\kappa R  < \delta_D(x_0) < R$ (later we will fix a point $z \in \partial D$
and restrict further that $x_0 \in B(z, R) \cap D$)
and set $\eps_1:=  \frac{\kappa R}{24}$. Define
$
 r(x,y): = \delta_D(x) \vee \delta_D(y)\vee |x-y|$  for $x,y\in D, $
and
\begin{equation} \label{d:gz1}
\BB(x,y):=
\begin{cases}
\left\{ A \in D:\, \delta_D(A) > \frac{\kappa}{2}r(x,y), \,
|x-A|\vee
|y-A| < 5 r(x,y)  \right\}& \text{ if } r(x,y) <\eps_1 \\
\{x_0 \}& \text{ if } r(x,y) \ge \eps_1 .
\end{cases}
\end{equation}

Set
$
C_2:=
C_1( 2) {\Phi(\frac{\delta_D(x_0)}{2})}{(\frac{\delta_D(x_0)}{2})^{-d}},
$
where $C_1(2)$ is the constant from Theorem \ref{t:greenh}.
By Lemma \ref{dec} we see that $G_D(x, x_0)\leq C_2$ for  $x\in
(D \cap B(x_0, 2)) \setminus B(x_0,
\frac{\delta_D(x_0)}{2})$.

Now we define
\begin{equation}\label{d:gz0}
g(x):=  G_D(x, x_0) \wedge C_2.
\end{equation}
We note that for $y \in D \cap B(x_0, 2)$ with $\delta_D(y) \le 6 \eps_1$, we have
$g(y )= G_D(y, x_0),$
since $6\eps_1<\frac{\delta_D(x_0)}{4}$ and thus
$|y-x_0| \ge \delta_D(x_0) - 6 \eps_1 \ge \frac{\delta_D(x_0)}{2}$.

The following lemma follows from
\cite[Theorem 2.10]{KSV10}.

\begin{lemma}\label{G:g2}
(1) There exists $c=c(\kappa, R, \Psi, \gamma)>1$
such that for every $x \in D\cap B(x_0, 2)$ satisfying
$\delta_D(x)\ge 2^{-6} \kappa^3\eps_1$ we have $c^{-1}  \le g(x) \le c\,.$

\noindent
(2) There exists $c=c(\kappa, R, \Psi, \gamma)>0$
such that for every $x, y \in D\cap B(x_0, 2)$,
$c^{-1} g(A_1) \le g(A_2) \le c g(A_1)$
for every $ A_1, A_2 \in \BB(x,y).$
\end{lemma}
\medskip

With these preparations, we can prove the following two-sided estimates on the Green functions of
bounded $\kappa$-fat open sets (without loss of generality, assuming that diam$(D)\le 1$),
which extend \cite[Theorem 1.2]{KSV2}. As we mentioned in \cite[Theorem 1.2]{KSV2},
with
\eqref{e:H1n}, Lemma \ref{dec}--\ref{G:g3}(1),
\cite[Theorems 2.10 and  2.18]{KSV10}
at hand,
one can easily adapt the arguments of \cite[Proposition 6]{B3}. Since these are
more or less standard now, we omit the details. (See also the proof of  \cite[Theorem 6.4]{KM}.)

\begin{thm}\label{t:Gest}
 If $D$ is a bounded $\kappa$-fat open set with
localization radius $R$,
then there exists a constant
$c=c($diam$(D), R,
\kappa, \Psi, \gamma)>1$
such that for every $x, y \in D$ and $A \in \BB(x,y)$,
\begin{equation}\label{e:Gest}
c^{-1}\,\frac{g(x) g(y)\Phi(|x-y|)}{g(A)^2|x-y|^{d}}\,\le\,
G_D(x,y) \,\le\, c\,\frac{g(x) g(y)\Phi(|x-y|)}{g(A)^2|x-y|^{d}},
\end{equation}
where $g$ and $\BB(x,y)$ are defined by \eqref{d:gz0} and
\eqref{d:gz1} respectively.
\end{thm}

In fact, it is the next result, which covers unbounded $\kappa$-fat open sets,
that we will use in the following section.

\begin{thm}\label{t:green}
 Suppose that  $D$ is a
 $\kappa$-fat open set with
localization radius $R$
and $z\in \partial D$.
Assume that $x_0 \in B(z, R) \cap D$ with
$\kappa R  < \delta_D(x_0) < R$.
There exists $C_3=C_3(\Psi, \gamma, R, \kappa)>1$
such that for all $x,y\in B(z, 2^{-7}\kappa^2 R) \cap D$ and $A \in \BB(x,y)$
it holds that
\begin{align}\label{e:green2}
 C_3^{-1} \frac{g(x) g(y)\Phi(|x-y|)}{g(A)^2|x-y|^{d}} \, \le\, G_D(x,y)
\,\le \, C_3
\frac{g(x) g(y)\Phi(|x-y|)}{g(A)^2|x-y|^{d}}\, .
\end{align}
\end{thm}

\pf
By Theorem \ref{t:Gest}, we only need to prove the
theorem for unbounded $D$.

First, note that $g(\cdot)=G_D(\cdot, x_0)$ on  $B(z, \kappa R/4)$.
Using Proposition \ref{p:localization}, we choose a bounded $\kappa$-fat open set $D_1$ with
localization radius $R$
such that $D\cap B(z, 1) \subset D_1 \subset D\cap B(z, 2)$.
First note that, since   $ \delta_{D_1} (x) \vee \delta_{D_1} (y)=\delta_D(x) \vee \delta_D(y)
\le 2^{-7}\kappa^2 R$,
by the boundary Harnack principle in \cite[Theorem 2.18(ii)]{KSV10},
\begin{align}\label{e:ggD}
\frac{(G_{D_1}(x, x_0)\wedge C_2) (G_{D_1}(y, x_0)\wedge C_2)}{(G_{D_1}(A, x_0)\wedge C_2)^2} = \frac{G_{D_1}(x, x_0)
G_{D_1}(y, x_0)}{G_{D_1}(A, x_0)^2}
\asymp \frac{g(x) g(y)}{g(A)^2}.
\end{align}
It follows from \eqref{e:ggD} and Theorem \ref{t:Gest} that
$$
G_D(x,y) \ge G_{D_1}(x,y) \ge    c_1 \frac{g(x) g(y)\Phi(|x-y|)}{g(A)^2|x-y|^{d}}.
$$
On the other hand, by the strong Markov property,  the symmetry of $G_D$,  \eqref{e:ggD}
and Theorem \ref{t:Gest},
\begin{align*}
G_D(x,y)=G_{D_1}(x,y)+\E_x[G_D(y, X_{\tau_{D_1}})]
\le c_2 \frac{g(x) g(y)\Phi(|x-y|)}{g(A)^2|x-y|^{d}} + \E_x[G_D(y, X_{\tau_{D_1}})].
\end{align*}
Thus, since $\inf_{a \le 1} \Phi(a)a^{-d} >0$,  it suffices to show the first inequality below:
\begin{align}\label{e:green3}
\E_x[G_D(y, X_{\tau_{D_1}})]\,\le \,c_3
\frac{g(x) g(y)}{g(A)^2}\le  c_4
\frac{g(x) g(y)\Phi(|x-y|)}{g(A)^2|x-y|^{d}} \, .
\end{align}

Let $\eta_0:=2^{-2}\kappa R$ and $\eta_1:=2^{-3}\kappa^2 R$.
Since
$$
|z-A| \le |x-z|+|A-x| \le 2^{-7}\kappa^2 R + 5 r(x,y) \le 11  \cdot 2^{-7}\kappa^2 R < 2^{-3}
\kappa^2 R=\eta_1,
$$
we have  $A \in D \cap B(z, \eta_1)$.
Thus by Lemma \ref{C:c_L},
\begin{equation}\label{e:gnew}
G_D(A, x_0) \le c_5 (G_{D}(A_{\eta_0}(z), x_0) \wedge G_{D}(A_{\eta_1}(z), x_0)).
\end{equation}
Applying
the boundary Harnack principle in \cite[Theorem 2.18(ii)]{KSV10}
to $G_D(\cdot, w)$ and $G_{D}(\cdot, x_0)$ and using \eqref{e:gnew}, we get
\begin{align}\label{e:green5}
\int_{D \setminus D_1} G_D(y, w) \P_x(X_{\tau_{D_1}} \in dw)
\le & c_6 \frac{G_{D} (y, x_0)}{G_{D}(A_{\eta_0}(z), x_0)} \int_{D \setminus D_1}
G_D(A_{\eta_0}(z), w)\P_x(X_{\tau_{D_1}} \in dw)\nn\\
 \le & c_5c_6 \frac{g (y)}{g(A)} \E_x[G_D(A_{\eta_0}(z), X_{\tau_{D_1}})].
\end{align}
Using
the boundary Harnack principle in \cite[Theorem 2.18(ii)]{KSV10}
 and \eqref{e:gnew} again,
\begin{align}\label{e:green7}
&\E_x[G_D(A_{\eta_0}(z), X_{\tau_{D_1}})] \le c_{7}
\frac{G_D(x, x_0)}{G_{D}(A_{\eta_1}(z), x_0)}
\E_{A_{\eta_1}(z)}[G_D(A_{\eta_0}(z), X_{\tau_{D_1}})] \nn\\
&\le c_{7}c_5  \frac{g(x)}{g(A)}
\E_{A_{\eta_1}(z)}[G_D(A_{\eta_0}(z),
 X_{\tau_{D_1}})]\, .
\end{align}
By the strong Markov property and Theorem \ref{t:greenh},
\begin{align}
&\E_{A_{\eta_1}(z)}[G_D(A_{\eta_0}(z), X_{\tau_{D_1}})] \le
\E_{A_{\eta_1}(z)}[G_D(A_{\eta_0}(z), X_{\tau_{D_1}})]+
G_{D_1}(A_{\eta_1}(z), A_{\eta_0}(z))\nn \\
&=G_{D}(A_{\eta_1}(z), A_{\eta_0}(z)) \le G(A_{\eta_1}(z), A_{\eta_0}(z))
\le c_8  \sup_{a \ge \eta_1} \frac{\Phi(a)} {a^d}:=c_9
< \infty\, . \label{e:green8}
\end{align}
Combining \eqref{e:green5}--\eqref{e:green8}, we have proved the first inequality in \eqref{e:green3}.
\qed

We remark in passing that one of the reasons we introduced the function $g$, instead of using only
the function $G_D(\cdot, x_0)$, is that the function $g$ satisfies the local scale invariant
Harnack inequality defined in Definition \ref{kernelHI} while the function $x \to G_D(x, x_0)$
does not.

%%%%%%%%%%%%%%%%%%%%%%%%%%%%%%%%%%%%%%%%%%%%%
%%%%%%%%%%%%%%%%%%%%%                      Martin kernel and estimates under H1           %%%%%%%%%%%%%%%%%%%%%%%%%%%%%%%%%%%%%%%%%%%%
%%%%%%%%%%%%%%%%%%%%%%%%%%%%%%%%%%%%%%%%%%%%%

\section{Martin kernel and estimates}\label{mke}

In this section we discuss Martin kernels and their estimates.
Let $D$ be an open set in  $\R^d$. Fix a point $x_0\in D$ and define
$$
M_D(x, y):=\frac{G_D(x, y)}{G_D(x_0, y)}, \qquad x, y\in D,~y\neq x_0.
$$
As the process $X^D$ satisfies Hypothesis (B) in \cite{KW}, $D$  has a Martin boundary $\partial_M D$
with respect to $X$ and $M_D(x ,\, \cdot\,)$ can be  continuously extended  to $\partial_M D$.
 The Martin kernel at $z\in \partial_M D$ (based at $x_0\in D$) is denoted by $M_D(x, z)$.

Recall that  a positive harmonic function $f$ for  $X^{D}$ is minimal if, whenever
$h$ is a positive harmonic function for $X^{D}$ with $h\le f$ on $D$,
one must have $f=ch$ for some constant $c$.

\begin{defn}\label{def:cMB}
(1) A point $z\in \partial_M D$ is called a finite Martin boundary point if there exists a bounded sequence
$\{w_n\}\subset D$ converging to $z$ in the Martin topology.

\noindent (2) A point $z\in \partial_M D$ is called an infinite Martin boundary point if every sequence
$\{w_n\}\subset D$ converging to $z$ in the Martin topology is unbounded (in the Euclidean topology).

\noindent (3) A point $z\in \partial_M D$ is called a minimal Martin boundary point
if $M_D(\cdot, z)$ is a minimal harmonic function for $X^D$.
Denote by $\partial_m D$ the minimal Martin boundary of $X^D$.
\end{defn}

In \cite[Section 3]{KSV10}, we showed that the finite part of the (minimal) Martin boundary of any $\kappa$-fat set $D$
coincides with the Euclidean boundary of $D$, see \cite[Theorem 3.13]{KSV10}.

Using Theorem \ref{t:green}, we get the following Martin kernel estimates.
Recall that  $g$ is defined by \eqref{d:gz0}.
\begin{thm}\label{t:mke}
Suppose that  $D$ is a $\kappa$-fat open set with
localization radius $R$
and $z \in \partial D$. Assume that  $x_0 \in B(z, R) \cap D$ with
$\kappa R  < \delta_D(x_0) < R$.
  There exists a constant
$C_4>0$ such that for all $x \in B(z, 2^{-7}\kappa^2 R) \cap D$,
$$
C_4^{-1}\frac{g(x)\Phi(|x-z|)}{g(A_{|x-z|}(z))^2|x-z|^{d}}\, \le\, M_D(x,z)
\,\le \, C_4
\frac{g(x)\Phi(|x-z|)}{g(A_{|x-z|}(z))^2|x-z|^{d}}\, .
$$
\end{thm}
\pf
Fix a point $x_1\in B(z, 2^{-7}\kappa^2 R) \cap D$. First we deal with the Martin kernel
$M^{x_1}_D$ based at $x_1$.
Since
$$
\delta_{D}(x) \vee \delta_{D}(y)\vee |x-y| \to  |x-z|   \quad \text{and} \quad
\delta_{D}(x_0) \vee \delta_{D}(y)\vee |x_1-y| \to |x_1-z|
$$
as $y \to z$, applying Theorem \ref{t:green} to $(x, y)$ and $(x_1, y)$ respectively, we get
$$
M^{x_1}_D(x,z)   \asymp
\frac{g(x) \Phi(|x-y|)|x_1-z|^{d}}
{g(x_1)\Phi(|x_1-z|)g(A_{|x-z|}(z))^2|x-z|^{d}} \asymp\frac{g(x)\Phi(|x-z|)}
{g(A_{|x-z|}(z))^2|x-z|^{d}}\, .
$$
The assertion of the theorem follows immediately from
the relationship  $M_D(x, \cdot)=M_D^{x_1}(x, \cdot)M_D(x_1, \cdot)$.
\qed

%%%%%%%%%%%%%%%%%%%%%%%%%%%%%%%%%%%%%%%%%%%%%%
%%%%%%%%%%%%%%%               Some results under H1 and H2             %%%%%%%%%%%%%%%%%%%%%%%%%%%%%%%%%%%%%%%%%%%%%%%%%
%%%%%%%%%%%%%%%%%%%%%%%%%%%%%%%%%%%%%%%%%%%%%%%

\section{Some results under {\bf (H1)} and {\bf (H2)}}\label{H1-H2}

In this section we assume that {\rm ({\bf H2})} also holds. We will
extend some known results and prove a new Green function estimate.
Our approach is heavily based on some recent results in \cite{BGR1, BGR2, CKi, CKS, CK, KSV10}.
Note that
({\bf H2}) implies transience of the process $X$
(since we are always assuming $d\ge 2$).

First note that if both ({\bf H1}) and ({\bf H2}) hold, there exist  $a_5, a_6>0$  such that
\begin{equation}\label{e:H2n}
a_5 \left(\frac{R}{r}\right)^{2(\delta_1 \wedge \delta_3)} \le \frac{\Psi(R)}{\Psi(r)} \le a_6 \left(\frac{R}{r}\right)^{2(\delta_2 \vee \delta_4)}, \quad a>0,\ 0<r<R<\infty\, ,
\end{equation}
cf.~\cite[(2.6)]{KSV8}.

It follows from \cite{BGR1, CK} that there exists a constant $C_5>1$ such that
\begin{align}\label{e:sbmG}
C_5^{-1}\frac{\Phi(|x|)}{|x|^d} \le G(x)\le C_5\frac{\Phi(|x|)}{|x|^d}, \quad \text{for all } x \in \R^d
\end{align}
and
\begin{align}\label{e:sbmJ}
C_5^{-1} \frac{1}{|x|^d\Phi(|x|)}\le J_X(x)\le C_5 \frac{1}{|x|^d\Phi(|x|)}, \quad \text{for all } x \in \R^d.
\end{align}

The next result is proved in a more general setting in \cite[Section 3]{CKi}.
In fact, one can also follow the proofs
in \cite[Section 3]{CKS} to see that
all the arguments of  \cite[Section 3]{CKS} with $T=\infty$ go through
using the (global) parabolic
Harnack inequality in \cite{CK},  \eqref{e:j-decay}, \eqref{e:psi1},
\eqref{e:H2n}, \eqref{e:sbmJ} and the semigroup property.
Thus by following the arguments in \cite[Section 3]{CKS} line by line,
one can also prove the next proposition.
We omit the details.

\begin{prop}\label{step31}
Suppose that  {\rm ({\bf H2})} holds.
Let $a$ be a positive constant.
There exists $c=c(a, \Psi, \gamma)>0$ such that for any open set $D$,
$
p_D(t, x, y)\ge  c  ((\Phi^{-1}(t))^{-d} \wedge  { t}{J_X(x-y)})
$
for every $(t, x, y)\in
(0, \infty)\times D\times D$ with
$ \delta_D(x) \wedge \delta_D (y) \ge a \Phi^{-1}(t)$.
\end{prop}

Using \eqref{e:sbmJ} and Proposition \ref{step31}, the proof of the next
lemma is almost identical to that of \cite[Lemma 2.14]{KSV10}.
We omit the details.

\begin{lemma}\label{G:g4}
Suppose that  {\rm ({\bf H2})} holds.
For every $L >0$ and every open set $D$, there exists $c=c(L, \Psi, \gamma)>0$ such
that for every $|x-y| \le L ( \delta_D(x) \wedge \delta_D(y)) $,
\begin{equation}\label{e:g4}
G_D(x,y) \ge c \frac{\Phi(|x-y|)}{|x-y|^d}.
\end{equation}
\end{lemma}

We now recall the following (global) scale invariant boundary Harnack inequality from \cite{KSV10} that will
be used in Section \ref{mti}.

\begin{thm}[{\cite[Theorem 2.18(i)(ii) and Remark 2.21]{KSV10}}]\label{t:bhp-all}
Suppose that {\rm ({\bf H2})} holds.
There exists $c= c(\Psi, \gamma)>0$ such that the following hold for all $r>0$.
\begin{itemize}
    \item[(i)]
For every $z_0 \in \R^d$, every open set $U \subset B(z_0,r)$ and for any nonnegative function $u$ in $\RR^d$ which is regular harmonic in $U$ with respect to $X$ and vanishes a.e.~in $U^c \cap B(z_0, r)$ it holds that
\begin{eqnarray*}
c^{-1}\E_x[\tau_{U}] \int_{B(z_0, r/2)^c}  j(|y-z_0|)  u(y)dy
\,\le \,u(x)\,\le\,
c \,\E_x[\tau_{U}] \int_{B(z_0, r/2)^c} j(|y-z_0|)  u(y)dy
\end{eqnarray*}
for every  $ x \in U \cap B(z_0, r/2)$.
   \item[(ii)]
For every $z_0 \in \R^d$, every open set $D\subset \R^d$ and any nonnegative functions $u, v$ in
$\R^d$ which are regular harmonic in $D\cap B(z_0, r)$ with respect to $X$ and
vanish a.e. in $D^c \cap B(z_0, r)$, we have
        $$
        \frac{u(x)}{v(x)}\,\le c^4\,\frac{u(y)}{v(y)}\, , \qquad \mbox{ for all } x, y\in D\cap B(z_0, r/2).
        $$
        \end{itemize}
\end{thm}

The next theorem  is a consequence of
\cite[Theorem 4.1 and Corollary 4.5]{BGR2}
and  Theorem \ref{t:bhp-all}(i).

\begin{thm}\label{L:2}
Suppose that  {\rm ({\bf H2})} holds and $\gamma=1$ in \eqref{e:psi1}.
There exists $c=c(\Psi)>0$  such that for every
$C^{1,1}$ with characteristics $(R, \Lambda)$,  $r \in (0, R]$,
$z\in \partial D$ and any nonnegative function $u$ in $\R^d$
which is harmonic in $D \cap B(z, r)$ with respect to $X$ and
vanishes continuously on $ D^c \cap B(z, r)$, we have
\begin{equation}\label{e:bhp_m}
\frac{u(x)}{u(y)}\,\le c\,
\sqrt{\frac{\Phi(\delta_{D}(x)
)}{\Phi(\delta_{D}(y)
)}}\ ,
\qquad \hbox{for every } x, y\in  D \cap B(z, \tfrac{r}{2}).
\end{equation}
\end{thm}

%%%%%%%%%%%%%%%%%%%%%%%%%%%%%%%%%%%%%%%%%
%%%%%%%%%%%%%%%%                Quasi-additivity of capacity             %%%%%%%%%%%%%%%%%%%%%%%%%%%%%%%%%%%%%%%%%%%%
%%%%%%%%%%%%%%%%%%%%%%%%%%%%%%%%%%%%%%%%

\section{Quasi-additivity of capacity}\label{qac}
In this section we will prove two types of results -- local and global.
For the global case, we  assume that {\bf (H2)} holds.
We always state the condition {\bf (H2)} explicitly
when we deal with the global case.

Let $\cp$ denote the capacity with respect to $X$ and $\cp_D$ the capacity with respect to the
killed process $X^D$. The goal of this section is to prove that $\cp_D$ is (locally) quasi-additive
with respect to Whitney decompositions of $D$.

We first revisit \cite[Section 5.4.1]{SV} and extend \cite[Proposition 5.55]{SV}.
Recall that  $G(x, y)$ (and $G_D(x,y)$) is the Green function of $X$ (and $X^D$, respectively), and let $G\mu(x)=\int G(x, y)\mu(dy)$ and
$G_D\mu(x)=\int_D G(x, y)\mu(dy)$.

For any compact subset $K$ of $\R^d$, let ${\cal P}_K$ be the set of
probability measures supported by $K$. Define
$$
e(K):=\inf_{\mu\in {\cal P}_K} \int G\mu(x)\, \mu(dx)\, .
$$
Since the kernel $G$ satisfies the Maria-Frostman maximum principle saying that
$\sup_{\R^d}G\mu=\sup_{\rm{Supp}(\mu)}G\mu$
(see, for example, Theorem 5.2.2~in \cite{Chu}), it follows from (\cite{Fug1}, page 159) that
for any compact subset $K$ of $\R^d$,
\begin{equation}\label{fuglede}
{\rm Cap}(K)=\frac{1}{\inf_{\mu \in{\cal P}_K}\sup_{x\in
{\rm Supp}(\mu)}G\mu(x)}=\frac{1}{e(K)}\, .
\end{equation}
Furthermore, the infimum is attained at the capacitary measure $\mu_K$.

Using \eqref{fuglede}, the proof of the next lemma is same as that of  \cite[Lemma 5.54]{SV}.
\begin{lemma}\label{capineq} Let $K$ be a compact subset of $\R^d$.
For any probability measure
$\mu$ on $K$, it holds that
\begin{equation}\label{hi honey}
\inf_{x\in\rm{Supp}(\mu)}G\mu(x) \le e(K)\le \sup_{x\in\rm{Supp}
(\mu)}G\mu(x) \, .
\end{equation}
\end{lemma}

\begin{prop}\label{capest}
There exist positive constants $c_1<c_2$
such that
\begin{equation}\label{capest formula}
c_1 \frac{r^{d}}{\Phi(r)}\le {\rm Cap}(\overline{B(0,r)}) \le  c_2 \frac{r^{d}}{\Phi(r)}
\quad \text{ for every } r \in (0, 1].
\end{equation}
Furthermore, if we assume {\rm ({\bf H2})}, \eqref{capest formula} holds for all $r>0$.
\end{prop}
\pf
We first consider the first claim.
By \eqref{fuglede} it suffices to show that
\begin{equation}\label{capest formula1}
c_1\frac{\Phi(r)}{r^{d}} \le e(\overline{B(0,r)}) \le  c_2 \frac{\Phi(r)} {r^{d}} \quad \text{ for every } r \in (0, 1].
\end{equation}

For every $w \in \overline{ B(0,r)}$, consider the intersection of
$B(0,r)$ and $B(w,r)$. This intersection contains the intersection
of $B(w,r)$ and the cone with vertex $w$ of aperture equal to $\pi/3$
pointing towards the origin.
Let $C(w)$ be the latter intersection.
Then by Theorem \ref{t:greenh},
$$
\int_{B(0,r)}G(w,y)dy  \ge
 \int_{C(w)}G(w,y)dy \ge c_1  \int_{B(w,r)}\frac{\Phi(|w-y|)}{|w-y|^{d}} dy=c_1 \int_{B(0,r)}\frac{\Phi(|y|)}{|y|^{d}} dy.
$$
Thus by {\bf (H1)},
\begin{equation}\label{e:gw1}
\inf_{w\in \overline{ B(0,r)}} \int_{B(0,r)}G(w,y)dy \ge c_1 \int_{B(0,r)} \frac{\Phi(|y|)}{|y|^{d}} dy
 \ge c_2\Phi(r).
\end{equation}

On the other hand, Theorem \ref{t:greenh} and  {\bf (H1)}
 also give \begin{equation}\label{e:gw2}
\sup_{w\in\overline{ B(0,r)}} \int_{B(0,r)}G(w,y)dy
 \le
c_3\sup_{w\in\overline{ B(0,r)}} \int_{B(0,r)}
\frac{\Phi(|w-y|)}{|w-y|^{d}} dy  \le
c_3 \int_{B(0,2r)}
\frac{\Phi(|y|)}{|y|^{d}} dy  \le
c_4 \Phi(r).
\end{equation}
By applying Lemma \ref{capineq} with the normalized Lebesgue measure on
$B(0,r)$,  the proposition now follows from \eqref{hi honey},
\eqref{e:gw1} and  \eqref{e:gw2}.

If we assume ({\bf H2}), we use \eqref{e:H2n}  and \eqref{e:sbmG} instead of Theorem \ref{t:greenh} and  {\bf (H1)} respectively. Then we get \eqref{capest formula1} for all $r>0$ by the same argument.
\qed

\bigskip

For any open set $D\subset \R^d$ let $\SS(D)$  denote the collection of all excessive functions with respect to $X^D$
and let $\SS^c(D)$ be the family of positive  functions in $\SS(D)$ which are continuous in the extended sense.
Recall that
positive harmonic functions with respect to $X^D$ are in $\SS^c(D)$.
For any $u\in \SS(D)$
and $E\subset D$, the reduced function of $u$ relative to $E$ in $D$ is defined by
\begin{align}\label{e:RE}
R^E_u(x)=\inf\{v(x): v\in \SS(D) \mbox{ and } v\ge u \mbox{ on } E\},
\quad x\in \R^d.
\end{align}
The lower semi-continuous regularization $\wh R^E_u$ of $R^E_u$ is called the balayage
of $u$ relative to $E$ in $D$. Since the process $X^D$ has a continuous transition density
$p_D(t,x,y)$, its semigroup is strongly Feller. Thus it follows easily from
\cite[Proposition V.2.2]{BH} that the cone of excessive functions $\SS(D)$ is a
balayage space in the sense of  \cite{BH}.

Given $u\in \SS^c(D)$, define a kernel $k_u:D\times D\to [0,\infty]$ by
\begin{align}\label{e:ku}
k_u(x,y):=\frac{G_D(x,y)}{u(x)u(y)}\, ,\quad x,y\in D\,.
\end{align}
Note that $k_u(x,y)$ is jointly lower semi-continuous on $D \times D$ by Remark \ref{r:lsc} and the assumptions that $u$ is positive and continuous in the extended sense.
 For a measure $\lambda$ on $D$ let $\lambda_u(dy):=\lambda(dy)/u(y)$. Then
$$
k_u\lambda (x):=\int_D k_u(x,y)\, \lambda(dy)=\int_D \frac{G_D(x,y)}{u(x)u(y)}\,
\lambda(dy)=\frac{1}{u(x)}\int_D G_D(x,y)\, \frac{\lambda(dy)}{u(y)}=\frac{1}{u(x)}G_D \lambda_u(dy)\, .
$$
We define a capacity with respect to the kernel $k_u$ as follows:
\begin{eqnarray*}
\sC_u(E):=\inf\{\|\lambda\|: k_u\lambda \ge 1 \textrm{ on }E\}\, , \quad E\subset D\, ,
\end{eqnarray*}
where $\|\lambda \|$ denotes the total mass of the measure $\lambda$ on $D$. The following dual
representation of the capacity of compact sets can be found in \cite[Th\'eor\`eme 1.1]{Fug}:
\begin{equation}\label{e:fuglede-equality}
\sC_u(K)=\sup\{\mu(K):\, \mu(D\setminus K)=0, k_u\mu \le 1\ \textrm{on }D\}\, .
\end{equation}
For a compact set $K\subset D$, consider  the balayage $\wh{R}_u^K$. Being a potential,
$\wh{R}_u^K=G_D \lambda^{K, u}$ for a measure $\lambda^{K, u}$ supported in $K$.
Define the Green energy of $K$ (with respect to $u$) as
$$
\gamma_u(K):=\int_D \int_D G_D(x,y) \lambda^{K, u}(dx)\, \lambda^{K, u}(dy) =
\int_D G_D \lambda^{K, u}(x)\, \lambda^{K, u}(dx)=\EE(G_D \lambda^{K, u},G_D \lambda^{K, u})\, .
$$
As usual, this definition of energy is extended first to open
and then to Borel subsets of $D$.
The following proposition relates the energy $\gamma_u$ with the capacity $\sC_u$.
\begin{prop}\label{p:gcu}
For all Borel subset $E\subset D$ it holds that
$
\gamma_u(E)=\sC_u(E).
$
\end{prop}
\pf
Clearly, it suffices to prove
the proposition
 for compact subsets $K$ of $D$.
Note first that
\begin{eqnarray*}
\sup\{\mu(K):\, \mu(D\setminus K)=0, k_u\mu \le 1\ \textrm{on }D\}
&=&\sup\{\mu(K):\, \mu(D\setminus K)=0, G_D \mu_u \le u\ \textrm{on }D\}\\
&=&\sup\{\int_D u(y)\, \lambda(dy)\, ,\lambda(D\setminus K)=0, G_D\lambda \le u\ \textrm{on }D\}\, .
\end{eqnarray*}
Since $\lambda^{K, u}(D\setminus K)=0$ and $G_D\lambda^{K, u}=\wh{R}_u^K\le u$ on $D$,
we conclude from the above and \eqref{e:fuglede-equality} that
$$
\sC_u(K)\ge \int_K u(y)\, \lambda^{K, u}(dy)\ge \int_K G_D\lambda^{K, u}(y)\,
\lambda^{K, u}(dy)=\gamma_u(K)\, .
$$
Conversely, $G_D\lambda^{K, u}=\wh{R}_u^K=  u$ on $K$, hence by the definition of $\sC_u(K)$ we have
$$
\sC_u(K)\le \int_D u(y)\, \lambda^{K, u}(dy)=\int_K G_D\lambda^{K, u}(y)\,
\lambda^{K, u}(dy)=\gamma_u(K)\, .
$$
\qed

Note that in case $u\equiv 1$, $\gamma_1(E)=\sC_1(E)=\cp_D(E)$.

Let $\{Q_j\}_{j\ge 1}$ be a Whitney decomposition of $D$.
Recall that $x_j$ is the center of $Q_j$.
For each $Q_j$, let $Q_j^*$ denote the interior of the double of $Q_j$.
Then $\{Q_j, Q_j^*\}$ is a quasi-disjoint decomposition of $D$
in the sense of \cite[pp.~146-147]{AE}.
\begin{defn}\label{kernelHI}
(1)  A kernel $k:D\times D\to [0,+\infty]$ is said to satisfy the local Harnack property
with respect to $\{Q_j, Q_j^*\}$ if
\begin{align}\label{e:hpk}
k(x,y)\asymp k(x',y)\, , \textrm{ for all }x,x'\in Q_j \textrm{  and all }y\in D\setminus Q_j^*\, ,
\end{align}
for all cubes $Q_j$ whose diameter is less than $r_1$ for some $r_1>0$
(with constants independent of the cubes).

\noindent
(2)  A kernel $k:D\times D\to [0,+\infty]$ is said to satisfy the Harnack property
with respect to $\{Q_j, Q_j^*\}$ if
\eqref{e:hpk} holds
for all cubes $Q_j$.
\end{defn}

\begin{defn}\label{d:sihi}
(1) A function $u:D\to (0,\infty)$ is said to satisfy the local scale invariant Harnack inequality
with respect to $\{Q_j\}$
if for some $r_1>0$ there exists $c_1=c_1(r_1)>0$ such that
\begin{equation}\label{e:u-scale-inv-har-local}
\sup_{Q_j}u\le c_1 \inf_{Q_j}u\, ,\quad \textrm{for all } Q_j \text{ with }\mathrm{diam}(Q_j)<r_1\, .
\end{equation}

\noindent
(2) A function $u:D\to (0,\infty)$ is said to satisfy the scale invariant Harnack inequality
with respect to $\{Q_j\}$ if there exists $c_2>0$ such that
\begin{equation}\label{e:u-scale-inv-har}
\sup_{Q_j}u\le c_2 \inf_{Q_j}u\, ,\quad \textrm{for all } Q_j \, .
\end{equation}
\end{defn}

\begin{lemma}\label{l:sihi-khp}
If $u\in \SS^c(D)$ satisfies the (local) scale invariant Harnack inequality with respect
to $\{Q_j\}$, then the kernel $k_u$ satisfies the (local) Harnack property with respect to $\{Q_j, Q_j^*\}$.
\end{lemma}
\pf
Note that for all $Q_j$ with $\mathrm{diam}(Q_j)< r_1$ (respectively for all $Q_j$),
the function $G_D(\cdot, y)$ is regular harmonic in $Q_j$ for every $y\in D\setminus Q_j^*$.
Together with the assumption that $u$ satisfies the (local) scale invariant Harnack inequality,
this proves the claim.
\qed

Typical examples of positive continuous excessive functions $u$ that satisfy the local scale invariant Harnack inequality
are functions $u\equiv 1$ and $u=G_D(\cdot,x_0)\wedge c$ with $x_0\in D$  and $c>0$ fixed.
Similarly, if $X$  satisfies ({\bf H2}),
the scale invariant Harnack inequality holds for the same functions,
see \cite[Theorem 4.12]{CK}.

We record now two lemmas.
\begin{lemma}\label{l:cap-comparability}
(1)
There exists a constant $c=c(\Psi,\gamma, r_1)\in (0,1)$ such that
\begin{equation}\label{e:cap-comparability}
c\, \cp_D(Q_j)\le \cp(Q_j)\le \cp_D(Q_j)
\end{equation}for all Whitney cubes whose diameter is less than $r_1$.

\noindent
(2) Suppose that {\rm({\bf H2})} holds.
Then there exists a
constant $c=c(
\Psi,\gamma)\in (0,1)$ such that \eqref{e:cap-comparability} holds for all Whitney cubes.
\end{lemma}
\pf
(1) By
\eqref{e:fuglede-equality} and Proposition  \ref{p:gcu}
we have that for every compact set $K\subset D$,
\begin{eqnarray*}
\cp_D(K)=\sup\{\mu(K):\, \mathrm{supp}(\mu)\subset K, G_D\mu \le 1 \text{ on }D\}\, .
\end{eqnarray*}
If $\mathrm{supp}(\mu)\subset K$ and $G\mu\le 1$ on $\R^d$, then clearly $G_D \mu\le 1$ on $D$.
This implies that $\cp(K)\le \cp_D(K)$ for all compact subset $K\subset D$,
in particular for each Whitney cube $Q_j$.

If $x,y\in Q_j$, then $|x-y|\le \mathrm{diam}(Q_j)<r_1$ and $|x-y|\le \mathrm{diam}(Q_j)
\le \mathrm{dist}(Q_j,\partial D)\le \delta_D(x)\wedge \delta_D(y)$.
Thus by Lemma \ref{G:g3}(2) and Theorem \ref{t:greenh}
 there exists $c\in (0,1)$ such that for all $Q_j$ of
diameter less than $r_1$ it holds that
\begin{equation}\label{e:comp-green-fn}
G_D(x,y)\ge c\,  G(x,y)\, ,\qquad x,y \in Q_j\, .
\end{equation}
Let $\mu$ be the capacitary measure for $Q_j$ (with respect to $X^D$), i.e., $\mu(Q_j)
=\cp_D(Q_j)$. Then by \eqref{e:comp-green-fn} for every $x\in Q_j$ we have
$$
1\ge  G_D \mu(x)=\int_{Q_j}G_D(x,y)\, \mu(dy)\ge \int_{Q_j}cG(x,y)\, \mu(dy)=G(c\mu)(x)\, .
$$
By the maximum principle it follows that $G(c\mu)\le 1$ everywhere on $\R^d$. Hence,
 $\cp(Q_j)\ge (c\mu)(Q_j)=c\, \cp_D(Q_j)$.

\noindent
(2) The proof is analogous: We use Lemma \ref{G:g4} instead of Lemma \ref{G:g3}(2),
and \eqref{e:sbmJ} instead of Theorem \ref{t:greenh}.

\qed

\begin{lemma}\label{l:comp-gamma-cap}
(1) Suppose that $u\in \SS^c(D)$ is a function satisfying the local scale invariant
Harnack inequality
\eqref{e:u-scale-inv-har-local}. Then for every $Q_j$ of diameter less than $r_1$ and every
$E\subset Q_j$ it holds that
\begin{equation}\label{e:comp-gamma-cap}
\gamma_u(E)\asymp u(x_j)^2 \cp_D(E)\, .
\end{equation}

\noindent
(2) Suppose that {\rm({\bf H2})} holds
and let  $u\in \SS^c(D)$ be a function satisfying the scale invariant
Harnack inequality
\eqref{e:u-scale-inv-har}. Then \eqref{e:comp-gamma-cap} holds for every $Q_j$
and every $E\subset Q_j$.
\end{lemma}

\pf
(1) It suffices to prove \eqref{e:comp-gamma-cap} for compact subsets
$K\subset E\subset Q_j$ and for $Q_j$ of diameter less than $r_1$.
Since  $k$ satisfies the local Harnack property  for $\{Q_j, Q_j^*\}$,  we have
$u\asymp u(x_j)$ on $Q_j$. Hence $\wh{R}_u^{K}\asymp u(x_j)\wh{R}_1^{K}$,
implying $G_D\lambda^{K, u}\asymp u(x_j)G_D \lambda^{K, 1}$ (everywhere on $D$). Therefore,
\begin{eqnarray*}
\gamma_u(K)&=&\int_{K} G_D\lambda^{K, u}(x)\, \lambda^{K, u}(dx) \asymp u(x_j)
\int_{K} G_D\lambda^{K, 1}(x)\, \lambda^{K, u}(dx)\nn\\
&=&u(x_j)\int_{K} G_D\lambda^{K, u}(x)\, \lambda^{K, 1}(dx)\asymp u(x_j)^2
\int_{K}\lambda^{K, 1}(dx)=u(x_j)^2 \mathrm{Cap}_D(K)\, .
\end{eqnarray*}

\noindent
(2) The proof is analogous to the proof of
(1). \qed

\begin{defn}\label{d:comp_measure}
Let  $\{Q_j\}$ be  a Whitney decomposition of $D$.

\noindent
(1) A Borel measure $\sigma$ on $D$
is locally comparable to the capacity $\sC_u$ with respect to $\{Q_j\}$ at $z \in \partial D$
if there exists $r_1, c_1>0$ such that
\begin{eqnarray*}
& &\sigma(Q_j)\asymp \sC_u(Q_j),\quad \textrm{for all }Q_j  \textrm{ with } Q_j \cap B(z, r_1)
\not= \emptyset\, ,\\
& &\sigma(E)\le c_1 \sC_u(E),\quad \textrm{for all Borel }E \subset D \cap B(z, 2 r_1).
\end{eqnarray*}

\noindent
(2) A Borel measure $\sigma$ on $D$
is comparable to the capacity $\sC_u$ with respect to $\{Q_j\}$ if there exists $c_2>0$ such that
$\sigma(Q_j)\asymp \sC_u(Q_j)$ for all $Q_j$, and $\sigma(E)\le c_2 \sC_u(E)$ for all Borel $E\subset D$.
\end{defn}

In order to construct a comparable measure we need (local) Hardy's inequalities. Recall that
the local Hardy inequality at $z\in \partial D$ and the Hardy inequality are introduced in
Definition \ref{def:LHI} and Definition \ref{def:HI}
respectively.

Define
$$
\sigma_u(E):=\int_E u(x)^2\Psi(\delta_D(x)^{-1})\, dx\, ,\quad E\subset D\, .
$$

\begin{prop}\label{p:comparable_capacity}
(1)  If $({\cal E}, {\cal F}_D)$ satisfies the local Hardy inequality at $z \in \partial D$, then
for any Whitney decomposition $\{Q_j\}$ of $D$
and any $u\in \SS^c(D)$ satisfying the local scale invariant Harnack inequality for $\{Q_j\}$,
$\sigma_u$ is locally comparable to the capacity $\sC_u$ with respect to $\{Q_j\}$ at $z$.

\noindent
(2)
Suppose that {\rm({\bf H2})} holds
and that $({\cal E}, {\cal F}_D)$ satisfies the Hardy inequality.
Then for any Whitney decomposition $\{Q_j\}$ of $D$
and any $u\in \SS^c(D)$ satisfying the  scale invariant Harnack inequality for $\{Q_j\}$,
$\sigma_u$ is comparable to the capacity $\sC_u$ with respect to $\{Q_j\}$.
\end{prop}

\pf
(1) Fix $z\in\partial D$ and let $r_1=r_0/2$ where $r_0$ is the constant in the Definition \ref{def:LHI}.
Since  $u$ satisfies the local scale invariant Harnack inequality for $\{Q_j\}$,
we have $u\asymp u(x_j)$ whenever diameter of $Q_j$ is less than $r_1$.
By Lemma \ref{l:comp-gamma-cap}(1)
we have that $\gamma_u(Q_j)\asymp u(x_j)^2 \cp_D(Q_j)$
whenever diameter of $Q_j$ is less than $r_1$.
On the other hand,
\begin{align}\label{e:4.3e2}
\sigma_u(Q_j)=\int_{Q_j} u(x)^2\Psi(\delta_D(x)^{-1})\, dx \asymp u(x_j)^2
\Psi((\mathrm{diam}(Q_j))^{-1})\,
 |Q_j|\, .
\end{align}
Lemma \ref{l:cap-comparability}(1)
and Proposition \ref{capest} imply that
$$\cp_D(Q_j)
\asymp \mathrm{Cap}(Q_j)\asymp
\frac{(\mathrm{diam} (Q_j))^{d}}{\Phi(\mathrm{diam} (Q_j))}  \asymp
 \Psi((\mathrm{diam}(Q_j))^{-1})\, |Q_j|
\qquad
\text{for all }Q_j \text{ with }Q_j\cap B(z;r_1)\neq \emptyset\, .
$$
Thus, $\gamma_u(Q_j)\asymp \sigma_u(Q_j)$.

Using the local Hardy inequality,  for any
Borel subset $E\subset D \cap B(z, 2r_1)$ and compact $K\subset E$,
\begin{eqnarray*}
\gamma_u(E)&\ge & \gamma_u(K)=\EE(G_D \lambda^{K, u}, G_D \lambda^{K, u})
\ge c_1 \int_K (G_D\lambda^{K, u})(x)^2 \Psi(\delta_D(x)^{-1})\, dx \\
&=& c_1
\int_K u(x)^2 \Psi(\delta_D(x)^{-1})\, dx =c_1 \sigma_u(K)\, .
\end{eqnarray*}
This proves that $\gamma_u(E)\ge c_1  \sigma_u(E)$.

Part (2) is proved analogously.
\qed

Now we can repeat the argument in the proof of \cite[Theorem 7.1.3]{AE} and conclude the following.
\begin{prop}\label{p:quasi-additivity}
(1)  If $({\cal E}, {\cal F}_D)$ satisfies the local Hardy inequality at $z \in \partial D$,
then for any Whitney decomposition $\{Q_j\}$ of $D$ and any $u\in \SS^c(D)$ satisfying the local scale invariant Harnack inequality for $\{Q_j\}$,
 the Green energy $\gamma_u$ is locally quasi-additive with respect to $\{Q_j\}$ at $z$:
There exists a positive constant $r_1>0$ such that
$$
\gamma_u(E)\asymp \sum_{j\ge 1} \gamma_u(E\cap Q_j)\, \quad \textrm{for all Borel }E
\subset D \cap B(z, r_1).
$$
(2)
Suppose that {\rm({\bf H2})} holds.
If $({\cal E}, {\cal F}_D)$ satisfies the Hardy inequality,
then for any Whitney decomposition $\{Q_j\}$ of $D$
and any $u\in \SS^c(D)$ satisfying the  scale invariant Harnack inequality for $\{Q_j\}$,
the Green energy $\gamma_u$
is quasi-additive with respect to $\{Q_j\}$:
$$
\gamma_u(E)\asymp \sum_{j\ge 1} \gamma_u(E\cap Q_j)\, ,\quad \textrm{for all Borel }E\subset D \, .
$$
\end{prop}
\pf
(1) Choose $r_1$ to be the constant from Definition \ref{d:sihi} and
let $ E\subset D\cap B(z, r_1)$. If $E\cap Q_j\neq \emptyset$, then
$\mathrm{diam}(Q_j)<\mathrm{dist}(Q_j,\partial D)<r_1$. By subadditivity of
$\gamma_u$, we have that $\gamma_u(E)\le \sum_j \gamma_u(E\cap Q_j)$.
For the reverse inequality we may assume that $\gamma_u(E)<\infty$. Then there
exists a measure $\mu$ such that $k_u \mu\ge 1$ on $E$ and $\|\mu\|\le 2 \gamma_u(E)$.
For each $Q_j$ such that $E\cap Q_j\neq \emptyset$, we decompose the measure $\mu$ into
$
\mu_j:=\mu_{|Q_j^*}$ and $\mu_j':=\mu_{|D\setminus Q_j^*}$.
Then either (i) $k_u \mu_j \ge \frac12 $ on $E\cap Q_j$, or (ii) $k_u \mu'_j (x)\ge \frac12 $ for
some $x\in E\cap Q_j$. Let $J_1$ denotes the set of indices $j$ for which (i) holds, and $J_2$
those for which (ii) holds. For $j\in J_1$ we have $\gamma_u(E\cap Q_j)\le 2\|\mu_j\|$.
Since the number of overlaps of $\{Q_j^*\}$ is uniformly bounded, it follows that
\begin{equation}\label{e:sum-over-J1}
\sum_{j\in J_1}\gamma_u(E\cap Q_j)\le 2\sum_{j\in J_1} \|\mu_j\| \le c_1 \|\mu\| \le 2c_1 \gamma_u(E)\, .
\end{equation}

For $j\in J_2$, by
the local Harnack property of $k_u$ (Lemma \ref{l:sihi-khp}) we have
$$
k_u \mu(y)\ge k_u \mu'_j(y) =\int_{D\setminus Q_j^*}k_u(y,w)\, \mu'_j(dw)\ge c_2
\int_{D\setminus Q_j^*}k_u(x,w)\, \mu'_j(dw)=c_2 k_u \mu(x)\ge \frac12 c_2\, .
$$
Therefore,
$$
k_u \mu \ge \frac12 c_2=: c_3^{-1} \, \qquad \text{on }\bigcup_{j\in J_2}Q_j\, ,
$$
implying that $\gamma_u(\cup_{j\in J_2}Q_j) \le c_3 \|\mu\|\le 2c_3 \gamma_u(E)$.
Since by Proposition \ref{p:comparable_capacity} $\sigma_u$ is locally comparable
to $\gamma_u$, it follows from the $\sigma$-additivity of $\sigma_u$ that
\begin{eqnarray*}
\sum_{j\in J_2}\gamma_u(E\cap Q_j)&\le & \sum_{j\in J_2}\gamma_u( Q_j)\le c
\sum_{j\in J_2} \sigma_u(Q_j)\\
&\le & c_4 \sigma_u\Big(\bigcup_{j\in J_2} Q_j\Big)\le c_5
\gamma_u\Big(\bigcup_{j\in J_2} Q_j\Big) \le c_6 \gamma_u(E)\, .
\end{eqnarray*}
Together with \eqref{e:sum-over-J1} this finishes the proof.

Part
(2) is proved analogously.
\qed

In the remainder of the section we discuss sufficient geometric conditions which imply
the (local) Hardy inequality.

For $v\in {\cal F}_D$,
\begin{align}\label{e:efork}
\EE(v,v)=\int_D\int_D (v(x)-v(y))^2 J_X(x-y)\, dy\, dx
+ 2\int_D v(x)^2 \kappa_D(x)\, dx\, ,
\end{align}
where  $\kappa_D$ is given by
$
\kappa_D(x):=\int_{D^c}
J_X(x-y)\, dy.
$

For $x\in D$, let $z_x$ be the point on $\partial D$ such that $|z_x-x|=\delta_D(x)$.
We say $D$ satisfies the local exterior volume condition at $z\in \partial D$
if there exist $r_0,c>0$ such that for every $x \in B(z, r_0) \cap D$,
 $|D^c\cap B(z_x,\delta_D(x))| \ge c \delta_D(x)^{d}$.

\begin{prop}\label{l:hardy2} The local Hardy inequality holds at $z \in \partial D$ if
 $D$ is an open set satisfying the local exterior volume condition at $z$.  \end{prop}
 \pf
Let $x \in B(z, r_0) \cap D$. Note that
$
\kappa_D(x) \asymp \int_{D^c}|x-y|^{-d}\Psi(|x-y|^{-1})dy.
$
If $y\in D^c \cap B(z_x, \delta_D(x))$, then $|x-y|\le |x-z_x|+|z_x-y|\le 2\delta_D(x)$.
Hence $\Psi(|x-y|^{-1})\ge \Psi(2^{-1}\delta_D(x)^{-1})\ge c_3 \Psi(\delta_D(x)^{-1})$.
This implies that
\begin{eqnarray*}
&&\int_{D^c} |x-y|^{-d}\Psi(|x-y|^{-1})\, dy \ge \int_{D^c\cap B(z_x,\delta_D(x))}
|x-y|^{-d}\Psi(|x-y|^{-1})\, dy\\
&&\quad \ge  c_3 \delta_D(x)^{-d} \Psi(\delta_D(x)^{-1})\, |D^c\cap B(z_x,\delta_D(x))|\ge c_4 \Psi(\delta_D(x)^{-1})\, ,
\end{eqnarray*}
where in the last inequality we used the local exterior volume condition at $z$.
Thus by \eqref{e:efork} we have that
$$
\EE(v,v)\ge \int_{D} v(x)^2\kappa_D(x)\, dx
\ge \int_{B(z, r_0) \cap D} v(x)^2\kappa_D(x)\, dx
\ge c_1 \int_{B(z, r_0) \cap D} v(x)^2\Psi(\delta_D(x)^{-1})\, dx\, .
$$
\qed

 We say $D$ satisfies the exterior volume condition if there exist $c>0$ such that for
every $x \in  D$,
 $|D^c\cap B(z_x,\delta_D(x))| \ge c \delta_D(x)^{d}$, where $z_x$ is a point
on $\partial D$ such that $|z_x-x|=\delta_D(x)$.

\begin{prop}\label{l:hardy}
Suppose that {\rm({\bf H2})} holds
and that $D$ is
either unbounded $\kappa$-fat open set whose upper Assouad dimension is strictly less than
$d-2(\delta_2 \vee \delta_4)$ (see {\rm \cite[Definition 2.1]{DV}} for the definition)
or an open set satisfying the exterior volume condition.
There exists a constant $c>0$ such that
$$
\EE(v,v)\ge c \int_D v^2(x)\Psi(\delta_D(x)^{-1})\, dx\, , \quad v\in \FF_D\, .
$$
\end{prop}
\pf
If $D$ is an unbounded $\kappa$-fat open set whose upper Assouad dimension is strictly less than $d-
2(\delta_2 \vee \delta_4)$ then it follows from \cite[Theorem 4 and Proposition 9]{DV} that
$$
 \EE(v,v) \ge \int_D\int_D (v(x)-v(y))^2
 J_X(x-y)\, dy\, dx  \ge c \int_D v^2(x)\Psi(\delta_D(x)^{-1})\, dx.
$$
Since
$
\kappa_D(x) \asymp \int_{D^c}|x-y|^{-d}\Psi(|x-y|^{-1})\, dy\, ,
$
for all $x\in D$, if $D$ is
an open set satisfying the exterior volume condition, the proof is similar to that of the previous proposition.
So we omit the proof.
\qed

%%%%%%%%%%%%%%%%%%%%%%%%%%%%%%%%%%%%%%%%%%%%%
%%%%%%%%%%%%%%%%               Minimal thinness at a finite Martin boundary point             %%%%%%%%%%%%%%%%%%%%%%%%%%%%%%%%%%%%%%%%
%%%%%%%%%%%%%%%%%%%%%%%%%%%%%%%%%%%%%%%%%%%

\section{Minimal thinness at a finite Martin boundary point}\label{mtf}

We start this section by recalling the definition of minimal thinness and
proving a general result for
minimal thinness of a set at any minimal Martin boundary point.

\begin{defn}\label{def:cMthin}
Let $D$ be an open set in  $\R^d$.
A set $E\subset D$
is said to be minimally thin in $D$ at $z\in  \partial_m D$  with respect to $X$
if $\wh R^E_{M_D(\cdot, z)}\neq M_D(\cdot, z)$.
\end{defn}

For any $z\in \partial_M D$, let $X^{D, z}=(X^{D, z}_t, \P^z_x)$
denote the $M_D(\cdot, z)$-process, Doob's $h$-transform
of $X^D$ with $h(\cdot)=M_D(\cdot, z)$. The lifetime of $X^{D, z}$ will be denoted by $\zeta$. It is known
(see \cite{KW}) that
$\lim_{t\uparrow\zeta}X^{D, z}_t=z$  $\P^z_x$-a.\/s..
For $E\subset D$, let $T_E:=
\inf\{t>0: X^{D, z}_t\in E\}$. It is proved in \cite[Satz 2.6]{Fol} that a set $E\subset D$ is minimally thin
at $z\in \partial_m D$ if and only if there exists $x\in D$ such that $\P^z_x(T_E<\zeta)\neq 1$.

The following proposition gives two more equivalent conditions for minimal thinness.

\begin{prop}\label{p:thinness}
Let $D$ be an open set in  $\R^d$,
$A\subset D$ and $z\in \partial_m D$. The following are equivalent:

\noindent
(1)
$A$ is minimally thin in $D$  at $z$ with respect to $X$.

\noindent
(2)
There exists an excessive function $u=G_D\mu +M_D \nu$  such that
$
\liminf_{x\to z, x\in A} \frac{u(x)}{M_D(x,z)}>0\, .$

\noindent
(3) There exists a potential $u=G_D\mu$ such that
$
\liminf_{x\to z, x\in A} \frac{u(x)}{M_D(x,z)}=+\infty \, .
$
\end{prop}
\pf We sketch the proof following the proof of \cite[Theorem 9.2.6]{AG}. Clearly,
(3) implies
(2).

Assume that
(2) holds. Then there exists a Martin topology neighborhood $W$ of $z$ and $a>\nu(\{z\})$
such that $u\ge a M_D(\cdot, z)$ on $A\cap W$.
If $\wh{R}^{A\cap W}_{M_D(\cdot, z)} =M_D(\cdot, z)$, then $u\ge \wh{R}^{A\cap W}_u
\ge a M_D(\cdot,z)$ everywhere. Thus $u-aM_D(\cdot, z)$ is excessive, hence $u-a
M_D(\cdot, z)=G_D\mu +M_D \wt{\nu}$ for a (non-negative) measure $\wt{\nu}$
on $\partial D$. On the other hand, $u-aM_D(\cdot, z)=G_D\mu +M_D
\nu_{|\partial D\setminus \{z\}}+(\nu(\{z\})-a)M_D(\cdot, z)$.
This implies that $\wt{\nu}=\nu_{|\partial D\setminus \{z\}}+(\nu(\{z\})-a)\delta_z$
yielding $\wt{\nu}(\{z\})=\nu(\{z\})-a<0$, which is a contradiction. Hence
$\wh{R}^{A\cap W}_{M_D(\cdot, z)} \neq M_D(\cdot,z)$, i.e., $A$ is
minimally thin at $z$. Thus
(2) implies
(1).

Suppose that
(1) holds. By \cite[Lemma 2.7]{Fol}, there exists an open subset
$U\subset \R^d$ such that $A\subset U$, and $U$ is minimally thin in $D$  at $z$
with respect to $X$. By the analog of \cite[Theorem 9.2.5]{AG}, there is a decreasing
sequence $(W_n)_{n\ge 1}$ of Martin topology open neighborhoods of $z$ shrinking to
$z$ and such that $\wh{R}^{U\cap W_n}_{M_D(\cdot,z)}(x_0)\le 2^{-n}$. Let $u_1:=
\sum_{n=1}^{\infty} \wh{R}^{U\cap W_n}_{M_D(\cdot,z)}$. Then $u_1$ is a sum of
potentials, hence a potential itself since $u_1(x_0)<\infty$. Further, $\wh{R}^{U\cap W_n}_{M_D(\cdot,z)}={M_D(\cdot,z)}$ on the open set $U\cap W_n$. Therefore,
$u_1(x)/M_D(x,z)\to \infty$ as
$x\to z$, $x\in U$. Thus
(3) holds. \qed

Note that this proposition holds true regardless whether $z$ is
a finite or an infinite Martin boundary point.

In the sequel we assume that
$D$ is a $\kappa$-fat open set with
localization radius $R$
and that $z$ is a fixed
point  in $\partial D$.
Without loss of generality we assume $R<1/10$.
Recall that $x_0 \in D\cap B(z, R)$ satisfies $\kappa R  < \delta_{D} (x_0)  < R$.
Let $M_D$ be the Martin kernel of $D$ based at $x_0$.

\begin{lemma}\label{l:key-step-1}
There exists $C_6=C_6(\Psi, \gamma, R, \kappa)>0$ such that
for every $x,y \in B(z, 2^{-7}\kappa^2 R) \cap D$ with
$|x-y|\ge \frac34 |x-z|$,
\begin{eqnarray}\label{e:key-step-1}
\frac{G_D(x,y)}{M_D(x,z)}\le C_6 G_D(x_0,y).
\end{eqnarray}
\end{lemma}
\pf
Recall that $C_3$ and $C_4$
are the constants from Theorem \ref{t:green} and
Theorem \ref{t:mke} respectively, and that $g(\cdot)=G_D(\cdot, x_0)$ on $B(z, \kappa R/4)$. We have
\begin{align}
\frac{G_D(x,y)}{M_D(x,z)}\le C_3  C_4
\frac{\Phi(|x-y|)|x-z|^{d}}{\Phi(|x-z|)|x-y|^{d}}
\frac{ G_D(y, x_0)G_D(A_{|x-z|}(z), x_0)^2}{G_D(A, x_0)^2}, \label{e:fpk1}
\end{align}
where $A \in \BB(x,y)$ and $\BB(x,y)$ is defined by \eqref{d:gz1}.
Since $|x-y|\ge \frac34 |x-z|$, we have
$
r(x,y)=
\delta_{D}(x) \vee \delta_{D}(y)\vee |x-y| \ge \frac34 |x-z|.
$

By the definition of $ A_{\frac43 r(x,y)}(z)$,
$\delta_D(A_{\frac43 r(x,y)}(z)) \ge \kappa \frac43 r(x,y) > \kappa r(x,y) /2$. Moreover,
$$
|x- A_{\frac43 r(x,y)}(z)|\, \le\, | x-z|+|z -  A_{\frac43 r(x,y)}(z)|
\,\le\, \frac{4}{3} |x-y| + \frac43 r(x,y) < 3 r(x,y)
$$
and
$ |y- A_{\frac43 r(x,y)}(z)| \le |y-x| +|x- A_{\frac43 r(x,y)}(z)|<4r(x, y)
< 5 r(x,y). $
Hence $A_{ \frac43 r(x,y)}(z)\in \BB(x,y).$

By Lemma
\ref{C:c_L} and  Lemma \ref{G:g2}(2),
\begin{align}
  G_D(A_{|x-z|}(z), x_0)
 \le c_1 G_D(A_{\frac43 r(x,y)}(z), x_0 )
= c_2 G_D(A_{\frac43 r(x,y)}(z), x_0 ) \le c_3 G_D(A, x_0).\label{e:fpk2}
\end{align}
 Moreover, by \eqref{e:Gwd},
$
\Phi(|x-y|)|x-z|^{d}/(\Phi(|x-z|)|x-y|^{d}) \le c_4.
$
The assertion of the lemma now follows from this, \eqref{e:fpk1} and \eqref{e:fpk2}.
\qed

The following proposition is an analog of \cite[Proposition V. 4.15]{BH}.
For $E\subset D$, define
$
E_n=E\cap \{x\in D:\, 2^{-n-1}\le |x-z|<2^{-n}\}\, ,\quad n\ge 1\, .
$
\begin{prop}\label{p:minthin-criterion-1}
A set
$E \subset D$ is minimally thin in $D$  at $z$ with respect to $X$ if and only if
 $
 \sum_{n=1}^{\infty}R^{E_n}_{M_D(\cdot, z)}(x_0)<\infty.
 $
\end{prop}
\pf
Assume $\sum_{n=1}^{\infty}R^{E_n}_{M_D(\cdot, z)}(x_0)<\infty$. Then there exists
$n_0\in \N$ such that
$\sum_{n_0}^{\infty}R^{E_n}_{M_D(\cdot, z)}(x_0)<\frac12 M_D(x_0,z)\, .$
Let $B=B(z,2^{-n_0})$. Then $A:=B\cap E\subset \cup_{n_0}^{\infty} E_n$. Therefore,
$R^A_{M_D(\cdot, z)}(x_0)\le \sum_{n_0}^{\infty}R^{E_n}_{M_D(\cdot, z)}(x_0)<
\frac12 M_D(x_0,z)$, implying $\wh{R}^A_{M_D(\cdot, z)}<\frac12 M_D(x_0,z)$. Hence,
there exists an excessive function $u$ such that $u\ge M_D(\cdot, z)$ on $A$ and $u(x_0)<
\frac12 M_D(x_0,z)$. By the Riesz decomposition, $s=G_D\mu +M_D\nu$. Hence,
$M_D \nu(x_0)=\int_{\partial_M D} M_D(x_0,y)\, \nu(dy)< \frac12 M_D(x_0,z)$
implying that $\nu(\{z\})<\frac12$.  Therefore,
$$
\liminf_{x\to z, x\in A} \frac{u(x)}{M_D(x,z)}\ge 1 >\frac12 >\nu(\{z\})\, .
$$
But this means that $A$ is minimally thin in $D$  at $z$ with respect to $X$. Clearly,
$E$ is also minimally thin in $D$  at $z$ with respect to $X$.

Conversely, suppose that $E$ is minimally thin in $D$  at $z$ with respect to $X$.
By Proposition \ref{p:thinness}, there exists a potential $u$ such that $
\liminf_{x\to z, x\in E} \frac{u(x)}{M_D(x,z)}=+\infty\, .
$
Without loss of generality, we may assume that
$
u(x_0)\le (2C_6)^{-1}.$
There exists $n_1\in \N$ with
$2^{-n_1} \le 2^{-7}\kappa^2 R$
 such that $u(x)>M_D(x,z)$ for all $x\in E\cap B(z,2^{-n_1})$. Thus, $E\subset
\overline{B}(z,2^{-n_1})^c\cup
\{u>M_D(\cdot,z)\}$.
For $n\ge n_1$ define
$$
F_n=\{x\in D:\, 2^{-n-1}<|x-z|<2^{-n}, u(x)>M_D(x,z)\}\qquad \textrm{and}
\qquad F=\bigcup_{n_1}^{\infty}F_n\, .
$$
Let $x\in E_n$. Since $|x-z|\le 2^{-n_1}$, we have that $u(x)>M_D(x,z)$ and
thus $x\in F_n$. This shows that $E_n\subset F_n$, $n\ge n_1$. Therefore, it suffices to show that
$\sum_{n_1}^{\infty}R^{F_n}_{M_D(\cdot, z)}(x_0)<\infty$.
Since $u>M_D(\cdot, z)$ on $F$,
it follows that $R^{F}_{M_D(\cdot, z)}(x_0)\le u(x_0)$.

Let $i\in\{1,2,3\}$. For every $n\in \N$, let
$
U_n=F_{n_1+3n+i}.
$
Since $i\in \{1,2,3\}$ is arbitrary, it suffices to show that $\sum_{n=1}^{\infty}
R^{U_n}_{M_D(\cdot, z)}(x_0)<\infty$. Let $U=\bigcup_{n=1}^{\infty}U_n$.
Then $U\subset F$ and thus $R^{U}_{M_D(\cdot, z)}(x_0)\le u(x_0)$. Note that
since $U$ is open, $\wh{R}^U_{M_D(\cdot, z)}=R^U_{M_D(\cdot, z)}$
(see \cite[page 205]{BH}).  Since $u$ is a potential, the same holds for
$\wh{R}^U_{M_D(\cdot, z)}$, hence there exists a measure $\mu$ such
that $R^U_{M_D(\cdot, z)}=G_D\mu$. Moreover, since $R^U_{M_D(\cdot, z)}$
is harmonic on $\overline{U}^c$ (cf.\cite[III.2.5]{BH}), $\mu(U^c)=0$. Let
$\mu_n:=\mu_{|\overline{U}_n}$. Since $\overline{U}_n$ are pairwise disjoint,
$$
\mu=\sum_{n=1}^\infty\mu_n \qquad \textrm{and} \qquad G_D\mu=
\sum_{n=1}^{\infty}G_D \mu_n\, .
$$
Fix $n\in \N$ and consider $l\in \N$, $x\in U_n$, $y\in \overline{U}_l$. If $l>n$,  then
$$
|y-z|<2^{-n_1-3l-i}\le 2^{1-3(l-n)}|x-z|\le \frac14 |x-z|\, ,
$$
and hence $|x-y|\ge |x-z|-|y-z|\ge \frac34 |x-z|$. If $l<n$, then analogously, $|x-z|\le
\frac14 |y-z|$, hence $|x-y|\ge |y-z|-|x-z|\ge \frac34 |y-z|\ge \frac34 |x-z|$. Thus, in both cases,
$$
|x-y|\ge \frac34 |x-z|\ge \frac34 \delta_D(x) \, ,\qquad x\in U_n, y\in \overline{U_l}, l\neq n\, .
$$

Define $\mu_n'=\mu-\mu_n$ and let $x\in U_n$. We have
\begin{eqnarray}
G_D\mu_n'(x)=\int_D G_D(x,y)\, \mu_n'(dy)= M_D(x,z)\int_D \frac{G_D(x,y)}{M_D(x,z)}\,
\mu_n'(dy).
\label{e:key-step-0}\end{eqnarray}
By \eqref{e:key-step-1}  we have that
\begin{eqnarray*}
&&G_D\mu_n'(x)
= M_D(x,z)\int_D \frac{G_D(x,y)}{M_D(x,z)}\, \mu_n'(dy) \le C_6 M_D(x,z) \int_D  G_D(x_0,y)\, \mu_n'(dy) \nonumber\\
&&\le  C_6 M_D(x,z) G_D \mu(x_0)\le C_6 M_D(x,z) u(x_0) <\frac12 M_D(x,z)\, .\nonumber
\end{eqnarray*}
Since $G_D\mu_n+G_D\mu_n'=G_D\mu=R^U_{M_D(\cdot, z)}=M_D(\cdot, z)$ on $U$, it follows that $G_D\mu_n=M_D(\cdot, z)-G_D\mu_n'\ge M_D(\cdot, z)-\frac12 M_D(\cdot, z)=\frac12
M_D(\cdot, z)$ on $U_n$. This implies that $G_D\mu_n \ge \frac12 R^{U_n}_{M_D(\cdot, z)}$. Finally,
$$
\sum_{n=1}^{\infty}R^{U_n}_{M_D(\cdot, z)}(x_0)\le 2\sum_{n=1}^{\infty}G_D
\mu_n(x_0)=2G_D\mu(x_0)<\infty\, .
$$
\qed

By Theorem \ref{t:mke}, for large $n$,
$$
C_4^{-1}2^{nd} \frac{
g(x)\Phi(2^{-n})}{g(A_{2^{-n}}(z))^2} \le M_D(x,z) \le C_4
2^{(n+1)d} \frac{
g(x)\Phi(2^{-n})}{g(A_{2^{-n}}(z))^2} \, ,\quad x\in E_n\, ,
$$

This implies that
$$
C_4^{-1}2^{nd} \frac{\Phi(2^{-n})}{g(A_{2^{-n}}(z))^2} R^{E_n}_{g} \le R^{E_n}_{M_D(\cdot, z)} \le C_4
2^{(n+1)d}\frac{\Phi(2^{-n})}{g(A_{2^{-n}}(z))^2} R^{E_n}_{g}\, .
$$
In particular,
\begin{equation}\label{e:equivalence-1}
\sum_{n=1}^{\infty}R^{E_n}_{M_D(\cdot, z)}(x_0)<\infty \quad \textrm{ if and only if}\quad \sum_{n=1}^{\infty}2^{nd}
\frac{\Phi(2^{-n})}{g(A_{2^{-n}}(z))^2}R^{E_n}_{g} (x_0)<\infty\, .
\end{equation}

Note that $\wh{R}^{E_n}_{g}$ is a potential, hence there exists a measure $\lambda_n$
(supported by $\overline{E}_n$) such that
$\wh{R}^{E_n}_{g}=G_D \lambda_n$. Also, $\wh{R}^{E_n}_{g}= g=G_D(\cdot, x_0)$
on $E_n$ (except for a polar set, and at
least for large $n$), hence
\begin{eqnarray*}
\wh{R}^{E_n}_{g}(x_0)&=& G_D \lambda_n(x_0)=\int_{\overline{E}_n}G_D(x_0,y)\, \lambda_n(dy)
=\int_{\overline{E}_n} g(y)\, \lambda_n(dy)\\
&=&\int_{\overline{E}_n}\wh{R}^{E_n}_{g}(y)\, \lambda_n(dy)=\int_D \int_D G_D(x,y)\,
\lambda_n(dy)\, \lambda_n(dx)= \gamma_{g}(E_n)\, .
\end{eqnarray*}
We conclude from \eqref{e:equivalence-1} that
\begin{equation}\label{e:equivalence-1-2}
\sum_{n=1}^{\infty}R^{E_n}_{M_D(\cdot, z)}(x_0)<\infty \quad \textrm{ if and only if}
\quad \sum_{n=1}^{\infty}2^{nd}\frac{\Phi(2^{-n})}{g(A_{2^{-n}}(z))^2} \gamma_{g}(E_n)<\infty\, .
\end{equation}
Thus we have proved the following Wiener-type criterion for minimal thinness.

\begin{corollary}\label{c:minthin-criterion-1}
$E \subset D$ is minimally thin in $D$  at $z$
with respect to $X$ if and only if
$$
\sum_{n=1}^{\infty}2^{nd} \frac{\Phi(2^{-n})}{g(A_{2^{-n}}(z))^2}\, \gamma_{g}(E_n)<\infty.
$$
\end{corollary}

Now we prove a version of Aikawa's criterion for minimal thinness.
\begin{prop}\label{p:aikawa-thinness}
Let $E\subset D$.
Let  $\{Q_j\}$ be  a Whitney decomposition of $D$ and let $x_j$ denote the center of $Q_j$.
If $({\cal E}, {\cal F}_D)$ satisfies the local Hardy inequality with a localization constant $r_0$
 at $z \in \partial D$,
then $E$ is minimally thin in $D$  at $z$ with respect to $X$ if and only if
$$
\sum_{ j: Q_j \cap B(z, r_0/2) \not= \emptyset} \mathrm{dist}(z,Q_j)^{-d}
\frac{\Phi(\mathrm{dist}(z,Q_j))}{g(A_{\mathrm{dist}(z,Q_j)}(z))^2}\, g(x_j)^2
\mathrm{Cap}_D(E\cap Q_j)<\infty\, .
$$
\end{prop}
\pf Let $r_1:=r_0/2$.
Further, let $V_n=\{x\in \R^d: 2^{-n-1}\le |x-z | <2^{-n}\}$ so that $E_n=E\cap V_n$.
If $V_n\cap Q_j \neq \emptyset$, then $\mathrm{dist}(z,Q_j)\asymp 2^{-n}$.
Since $g$ satisfies the local scale invariant Harnack inequality,
by the local quasi-additivity of $\gamma_{g}$ at $z$ (Proposition \ref{p:quasi-additivity}(1)),
\begin{eqnarray*}
&&\sum_{n=1}^{\infty}2^{nd} \frac{\Phi(2^{-n})}{g(A_{2^{-n}}(z))^2}\gamma_{g}(E_n)\,
\asymp \, \sum_{n=1}^{\infty}2^{nd}\frac{\Phi(2^{-n})}{g(A_{2^{-n}}(z))^2}
\sum_{Q_j \cap B(z, r_1) \not= \emptyset} \gamma_{g}(E_n\cap Q_j)\\
&\asymp &\sum_{j: Q_j \cap B(z, r_1) \not= \emptyset} \sum_{n: V_n\cap Q_j\neq
\emptyset} \mathrm{dist}(z,Q_j)^{-d}\frac{\Phi(\mathrm{dist}(z,Q_j))}
{g(A_{\mathrm{dist}(z,Q_j)}(z))^2} \gamma_{g}(E_n\cap Q_j)\\
&=&\sum_{Q_j \cap B(z, r_1) \not= \emptyset} \mathrm{dist}(z,Q_j)^{-d}
\frac{\Phi(\mathrm{dist}(z,Q_j))}{g(A_{\mathrm{dist}(z,Q_j)}(z))^2}
\sum_{n: V_n\cap Q_j\neq \emptyset}\gamma_{g}(E_n\cap Q_j)\\
&\asymp &\sum_{Q_j \cap B(z, r_1) \not= \emptyset} \mathrm{dist}(z,Q_j)^{-d}
\frac{\Phi(\mathrm{dist}(z,Q_j))}{g(A_{\mathrm{dist}(z,Q_j)}(z))^2}\, \gamma_{g}(E\cap Q_j)\, .
\end{eqnarray*}
In the second line above we used the fact that $g(A_{2^{-n}}(z))$ and
$g(A_{\mathrm{dist}(z,Q_j)}(z))$
are comparable, which is a consequence of
\cite[Theorem 2.10]{KSV10}.
For the last line we argue as follows: One inequality is the subadditivity of
capacity. For the other note that there exists $N\in \N$ such that for every $Q_j$,
$\sum_{n, V_n\cap Q_j\neq \emptyset}1=\sum_n 1_{V_n\cap Q_j}\le N$. Hence,
$\sum_{n, V_n\cap Q_j\neq \emptyset}\gamma_{g}(E\cap V_n\cap Q_j)\le
\sum_{n, V_n\cap Q_j\neq \emptyset}
\gamma_{g}(E \cap Q_j)\le N \gamma_{g}(E \cap Q_j)$.

Finally, by Lemma \ref{l:comp-gamma-cap} we see that
$\gamma_{g}(E\cap Q_j)\asymp g(x_j)^2 \cp_D(E\cap Q_j)$
which completes the proof by
Corollary  \ref{c:minthin-criterion-1}.
\qed

The next result is an analog of \cite[Part II, Corollary 7.4.4]{AE}.

\begin{corollary}\label{c:aikawa-thinness}
Suppose that either (i) $D$ is a half space; or (ii) $D\subset \R^d$ is a $C^{1,1}$ open set and $\gamma=1$
in \eqref{e:psi1}.
Let $x_j$ denote the center of $Q_j$. Then $E$ is minimally thin in $D$  at $z\in \partial D$
 with respect to $X$ if and only if
$$
\sum_{j: \, Q_j \cap B(z, 1) \not= \emptyset} \mathrm{dist}(z,Q_j)^{-d}
\Phi(\mathrm{dist}(Q_j, \partial D)) \mathrm{Cap}_D(E\cap Q_j)<\infty\, .
$$
\end{corollary}
\pf
The function $g$ is harmonic in $D\cap B(z, 2r_1)$ where $r_1:=  \kappa R/4$.
Since $X$ satisfies ({\bf H1}), applying
\cite[Theorem 2.18(i)]{KSV10}, we get that for $Q_j \cap B(z, r_1/10) \not= \emptyset$,
\begin{align}\label{e:newq1}
g(A_{\mathrm{dist}(z,Q_j)}(z))
\asymp \E_{A_{\mathrm{dist}(z,Q_j)}(z)}[\tau_{D\cap B(z, 2r_1)}]
\int_{B(z, r_1)^c}j(|y-z|)g(y)dy
\end{align}
and
\begin{align}\label{e:newq2}
g(x_j) \asymp \E_{x_j}[\tau_{D\cap B(z, 2r_1)}]
\int_{B(z, r_1)^c}j(|y-z|)g(y)dy.
\end{align}

Suppose that $D$ is a $C^{1,1}$ open set and $\gamma=1$. Since $C^{1,1}$
 set satisfies the interior and exterior ball condition, by combining
\eqref{e:newq1} and  \eqref{e:newq2} with \cite[Theorem 4.1 and Corollary 4.5]{BGR2}, we get
$$
g(A_{\mathrm{dist}(z,Q_j)}(z)) \asymp
\Phi(\mathrm{dist}(z,Q_j))^{1/2} \quad \text{and}\quad
g(x_j) \asymp \Phi(\mathrm{dist}(Q_j, \partial D))^{1/2}.
$$

In case when $D$ is the half space $\H$, the two relations above are immediate consequence of
\cite[Proposition 2.4]{BGR2}, the boundary Harnack principle in \cite[Theorem 2.18(ii)]{KSV10}
and the fact that $V(x):=V(x_d)$, where $V$ is the renewal function of the ascending
ladder height process of the $d$-th
component
 of $X$, is harmonic in $\H$ with respect to $X$.
Thus the corollary follows immediately from Proposition \ref
{p:aikawa-thinness}.

\qed

\noindent
{\bf Proof of Theorem \ref{t:dahlberg}:} Let $r_1:=r_0/2$ and without loss of
generality we assume $r_0< \kappa R/4$.  Assume that $E$ is minimally thin in $D$
at $z$ with respect to $X$.
Recall that $g(x ):=  G_D(x, x_0) \wedge C_2$ and $g(\cdot)=G_D(\cdot, x_0)$ on $B(z,\kappa R/4)$.
Thus, by Proposition \ref{p:aikawa-thinness},
$$
\sum_{j: Q_j \cap B(z, r_1) \not= \emptyset} \mathrm{dist}(z,Q_j)^{-d}
\frac{\Phi(\mathrm{dist}(z,Q_j))}{G_D(A_{\mathrm{dist}(z,Q_j)}(z), x_0)^2}G_D(x_j, x_0)^2
\mathrm{Cap}_D(E\cap Q_j)<\infty\, .
$$
First note that for $Q_j \cap B(z, r_1) \not= \emptyset$,
$$
\mathrm{Cap}_D(E\cap Q_j)\ge c_1 \sigma(E\cap Q_j)=c_1\int_E \1_{Q_j}(x)
\Psi(\delta_D(x)^{-1})\, dx\, .
$$
Next,
$\mathrm{dist}(z,Q_j)\asymp |x-z|$ and $G_D(x_j, x_0) \asymp G_D(x, x_0)$ for $x\in Q_j$. Therefore,
\begin{eqnarray*}
&&\sum_{Q_j \cap B(z, r_1) \not= \emptyset} \mathrm{dist}(z,Q_j)^{-d}
\frac{\Phi(\mathrm{dist}(z,Q_j))}{G_D(A_{\mathrm{dist}(z,Q_j)}(z), x_0)^2} G_D(x_j, x_0)^2
\mathrm{Cap}_D(E\cap Q_j) \\
&\ge &c_2 \sum_{Q_j \cap B(z, r_1) \not= \emptyset} \int_E |x-z|^{-d}
\frac{\Phi(|x-z|)}{G_D(A_{|x-z|}(z), x_0)^2} G_D(x, x_0)^2 \1_{Q_j}(x)\Psi(\delta_D(x)^{-1})\, dx\\
&\asymp & \sum_{Q_j \cap B(z, r_1) \not= \emptyset} \int_E |x-z|^{-d}
\left( \frac{G_D(x, x_0)}{G_D(A_{|x-z|} (z), x_0 )}\right)^2 \frac{\Psi(\delta_D(x)^{-1})}
{\Psi(|x-z|^{-1})}  \1_{Q_j}(x) dx\\
&=& \int_E \sum_{Q_j \cap B(z, r_1) \not= \emptyset} |x-z|^{-d}
\left( \frac{G_D(x, x_0)}{G_D(A_{|x-z|} (z), x_0 )}\right)^2 \frac{\Psi(\delta_D(x)^{-1})}
{\Psi(|x-z|^{-1})} \1_{Q_j}(x)\, dx\\
& \ge & \int_{E\cap B(z,r_1)} \left( \frac{G_D(x, x_0)}{G_D(A_{|x-z|} (z), x_0 )}\right)^2
\frac{\Psi(\delta_D(x)^{-1})}{\Psi(|x-z|^{-1})} |x-z|^{-d}\, dx\, .
\end{eqnarray*}
Conversely, assume that $E$ is the union of a subfamily of Whitney cubes.
Then $E\cap Q_j$ is either empty or equal to $Q_j$. Since $\mathrm{Cap}_D(Q_j)
\asymp \sigma(Q_j)=\int_{Q_j} \Psi(\delta_D(x)^{-1})\, dx$ for $Q_j \cap B(z, r_1)
\not= \emptyset$ by Proposition \ref{p:quasi-additivity}(1), we can reverse the first
inequality in the display above to conclude that
\begin{eqnarray*}
&&
\sum_{j:Q_j \cap B(z, r_1) \not= \emptyset} \mathrm{dist}(z,Q_j)^{-d}
\frac{\Phi(\mathrm{dist}(z,Q_j))}{G_D(A_{\mathrm{dist}(z,Q_j)}(Q), x_0)^2}G_D(x_j, x_0)^2
\mathrm{Cap}_D(E\cap Q_j)\\
&\le& c_3 \int_{E\cap B(z,5r_1)}\left( \frac{G_D(x, x_0)}{G_D(A_{|x-z|} (z), x_0)}\right)^2
\frac{\Psi(\delta_D(x)^{-1})}{\Psi(|x-z|^{-1})}  |x-z|^{-d}\, dx\, .
\end{eqnarray*}
\qed

\noindent
{\bf Proof of Corollary \ref{c:dahlberg}:}
We have seen
from the proof of Corollary \ref{c:aikawa-thinness} that there exists $r>0$ such that
for $Q_j \cap B(z, r) \not= \emptyset$,
$$
g(A_{\mathrm{dist}(z,Q_j)}(z))=G_D(A_{\mathrm{dist}(z,Q_j)}(z), x_0) \asymp
\Phi(\mathrm{dist}(z,Q_j))^{1/2}
$$
and
$$
g(x_j)=G_D(x_j, x_0) \asymp \Phi(\mathrm{dist}(Q_j, \partial D))^{1/2}.
$$
Combining the two relations above with the proof of Theorem \ref{t:dahlberg}, we immediately
arrive at the conclusion of Corollary \ref{c:dahlberg}.
\qed

%%%%%%%%%%%%%%%%%%%%%%%%%%%%%%%%%%%%%%%%%%%%%%%
%%%%%%%%%%%%%%%%%%                   Minimal thinness at infinity           %%%%%%%%%%%%%%%%%%%%%%%%%%%%%%%%%%%%%
%%%%%%%%%%%%%%%%%%%%%%%%%%%%%%%%%%%%%%%%%%%%%

\section{Minimal thinness at infinity}\label{mti}

Throughout this section we assume that {\rm ({\bf H2})} holds and the constant
$\gamma$ in \eqref{e:psi1} is $1$.
Thus $X$ is a unimodal L\'evy process satisfying the global weak scaling conditions in \cite{BGR2, BGR3}.
We will establish results for minimal thinness at infinity.
Even though the results are analogous to those of the previous section,
their proofs contain non-trivial modifications. In particular we will use the recently
established boundary Harnack principles given in  Theorems \ref{t:bhp-all} and \ref{L:2}.
Thus we include all details except in the proof of Theorem \ref{t:mbatinfty}.

We first extend the main result in \cite{KSV9}.
Let $\kappa \in (0,1/2]$. We say that an open set $D$ in $\R^d$ is $\kappa$-fat at infinity
if there exists $R>0$ such that for every $r\in [R,\infty)$ there exists $A_r \in \R^d$ such
that $B(A_r, \kappa r)\subset D\cap \overline{B}(0,r)^c$ and $|A_r|< \kappa^{-1} r$.

\begin{thm}\label{t:mbatinfty}
The Martin boundary at infinity with respect to $X$ of any open set $D$
which is $\kappa$-fat at infinity consist of exactly one point $
\infty$. The point  $\infty$
is a minimal Martin boundary point.
\end{thm}
\pf
The theorem is proved in \cite{KSV9} when $X$ is a subordinate Brownian motion with L\'evy exponent
$\Psi(\xi)=\phi(|\xi|^2)$  where $\phi$ is
a complete Bernstein function satisfying ({\bf H1}) and ({\bf H2}).
The method in \cite{KSV9} is quite robust and can be applied to
unimodal L\'evy processes satisfying the global weak scaling conditions. In fact, since we have
\eqref{e:H2n},  \cite[Lemma 2.2]{KSV10}, \eqref{e:sbmG} and \eqref{e:sbmJ}
(instead of
\cite[(2.2), Lemma 2.2, (2.8) and (2.9)]{KSV9} respectively),
using Theorems \ref{t:bhp-all} and \ref{L:2} instead of \cite[Theorem 1.1]{KSV7}
and \cite[Theorem]{KSV8}, one can follow the proofs in  \cite[Section 3]{KSV9} line by line
and obtain the corresponding results in  \cite[Section 3]{KSV9}.
Once we get the corresponding results in  \cite[Section 3]{KSV9},
then all arguments and results in \cite[Section 4]{KSV9} stay the same so that the theorem holds.
We omit the details since
these would be a simple repetition of proofs in \cite{KSV9}.
\qed

Since half-space-like open sets
are $\kappa$-fat at infinity, the Martin boundary at
infinity with respect to $X$ of any half-space-like open set
consists of exactly one point $\infty$ and this point is a minimal Martin boundary point.

In the remainder of this section we assume that $D\subset \R^d$ is a half-space-like open set
and that $
\bH_{1}\subset D\subset \bH.
$
Let $x_0=(\wt 0, 5)$ and let $M_D$ be the Martin kernel of $D$ based at $x_0$.

Before we prove Proposition \ref{p:minthin-criterion-1-infty}, which is an analog of
Proposition \ref{p:minthin-criterion-1} at infinity, we establish an inequality involving
Green functions and Martin kernel at infinity.
We recall from \cite{BGR3,CKi} that
for the half space $\H$ we have the following estimates: There exists a constant $c\ge 1$ such that
\begin{eqnarray}
  G_{\H}(x,y) \asymp \frac{\Phi(|x-y|)}{|x-y|^d}
 \left(1\wedge \frac{\Phi(\delta_{\H}(x))}{\Phi(|x-y|)}\right)^{1/2}\left(1\wedge
\frac{\Phi(\delta_{\H}(y))}{\Phi(|x-y|)}\right)^{1/2}\, . \label{e:green-function-estimate}
\end{eqnarray}
This implies
\begin{equation}\label{e:martin-kernel-estimate-infty}
c^{-1}\Phi(\delta_{\H}(x))^{1/2}\le M_\H(x,\infty) \le c \Phi(\delta_{\H}(x))^{1/2}\, .
\end{equation}
Relation \eqref{e:green-function-estimate} also implies that for every $u, v \in \H_2$,
\begin{align}\label{eGh01}
G_{\H_1}(u,v) \ge
 c_1\frac{\Phi(|u-v|)}{|u-v|^d} \left(1\wedge
\frac{\Phi(\delta_{\H}(u))}{\Phi(|u-v|)}\right)^{1/2}\left(1\wedge
\frac{\Phi(\delta_{\H}(v))}{\Phi(|u-v|)}\right)^{1/2} \ge c_2 G_{\H}(u,v).
\end{align}
Moreover, if $|x| \ge 10$, then $|x-x_0| \asymp |x| \ge \delta_{\H}(x)$. Thus for  $|x| \ge 10$,
\begin{align}\label{e:g-M-estimate-23}
G_{\H}(x, x_0) \asymp \frac{\Phi(|x|)}{|x|^d} \left(1\wedge
\frac{\Phi(\delta_{\H}(x))}{\Phi(|x|)}\right)^{1/2}\left(1\wedge \frac{1}{\Phi(|x|)}\right)^{1/2}
\asymp
\frac{\Phi(\delta_{\H}(x))^{1/2} }{|x|^d}.
\end{align}

\begin{lemma}\label{l:Gfeinfinty}
There exists $C_7>0$ such that for every $x,y \in B(0, 10)^c \cap D$ with
$|x-y|\ge \frac34 |y|$,
\begin{align}\label{e:Gfeinfinty}
{G_D(x,y)} \le C_7{G_D(x_0,y)}{M_D(x, \infty)}.
\end{align}
\end{lemma}
\pf
We claim that for every $w \in B(0, 10(|x| \vee |y|))^c \cap D$ with $\delta_D(w) \ge 3$,
\begin{align}\label{e:Gfeinfinty1}
\frac{G_D(x,y)}{G_D(x,w)} \le c_1 \frac{G_D(x_0,y)}{G_D(x_0,w)}.
\end{align}
By letting $w \to \infty$ with $\delta_D(w) \ge 3$, this implies \eqref{e:Gfeinfinty} immediately.

We prove \eqref{e:Gfeinfinty1} through 3 steps.

\noindent {\it Step 1.} We first prove \eqref{e:Gfeinfinty1} for $\H$.
Since
$|x-w| \asymp |x_0-w|$ and $\delta_{\H}(x) \le |x| \le |w| -|x| \le |x-w|$,
by \eqref{e:green-function-estimate},
\begin{align*}
&\frac{G_{\H}(x,y)G_{\H}(x_0,w)}{G_{\H}(x,w)} \nn\\
&\le c_2\frac{\Phi(|x-y|)}{|x-y|^d} \left(\frac{\Phi(\delta_{\H}(x))}
{\Phi(|x-y|)}\right)^{1/2}\left( \frac{\Phi(\delta_{\H}(y))}
{\Phi(|x-y|)}\right)^{1/2}\left(\frac{\Phi(\delta_{\H}(x_0))}{\Phi(|x_0-w|)}\right)^{1/2}
\left(\frac{\Phi(\delta_{\H}(x))}{\Phi(|x-w|)}\right)^{-1/2}\nn\\
&\le c_3\frac{\Phi(\delta_{\H}(x))^{1/2}\Phi(\delta_{\H}(y))^{1/2}}
{|x-y|^{d}\, \Phi(\delta_{\H}(x))^{1/2}}=
c_3\frac{\Phi(\delta_{\H}(y))^{1/2}}{|x-y|^{d}}.
\end{align*}
Thus, by our assumption $|x-y|\ge \frac34 |y|$ and \eqref{e:g-M-estimate-23},
$$
\frac{G_{\H}(x,y)G_{\H}(x_0,w)}{G_{\H}(x,w)} \le
c_4 |y|^{-d}
\Phi(\delta_{\H}(y))^{1/2} \le c_5 G_{\H}(x_0,y).
$$
We have proved \eqref{e:Gfeinfinty1} for $\H$.

\noindent {\it Step 2.} We assume $\delta_D(x) \wedge \delta_D(y) \ge 3$.
Then using the monotonicity of Green functions \eqref{eGh01} and Step 1, we have
$$
\frac{G_D(x,y)}{G_D(x,w)} \le \frac{G_{\H}(x,y)}{G_{\H_1}(x,w)} \le c_6
\frac{G_{\H}(x,y)}{G_{\H}(x,w)}
\le c_{7} \frac{G_{\H}(x_0,y)}{G_{\H}(x_0,w)}
\le c_{7} c_6^{-1} \frac{G_{\H_1}(x_0,y)}{G_{\H}(x_0,w)}\le  c_{7} c_6^{-1}
\frac{G_{D}(x_0,y)}{G_{D}(x_0,w)} .
$$
\noindent {\it Step 3.}
Let
\begin{align}
x_1:=
\begin{cases}
x & \text{if }  \delta_D(x)>3\\
(\wt x, 3) & \text{if }  \delta_D(x) \le 3
\end{cases}
\qquad \text{ and } \qquad
y_1:=
\begin{cases}
y & \text{if }  \delta_D(y)>3\\
(\wt y, 3) & \text{if }  \delta_D(y) \le 3.
\end{cases}
\end{align}
We use Theorem \ref{t:bhp-all}  when $\delta_D(x) \le 3$ and get
\begin{align}\label{e:gge13}
G_D(x,w)
=\frac{G_D(x,w)}{G_D(x, x_0)}G_D(x_0, x) \ge c_{11} \frac{G_D(x_1,w)}{G_D(x_1, x_0)}G_D(x_0, x).
\end{align}
Since $|x-y|\ge \frac34 |y| \ge\frac{15}{2} $,  we use Theorem \ref{t:bhp-all}  when $\delta_D(x) \le 3$
and get
\begin{align}\label{e:gge14}
G_D(x,y)
=\frac{G_D(x,y)}{G_D(x, x_0)}G_D(x_0, x) \le c_{12} \frac{G_D(x_1,y)}{G_D(x_1, x_0)}G_D(x_0, x).
\end{align}
If
 $\delta_D(x) \le 3$ then
$|y-x_1|
\ge |y-x|-|x-x_1|  \ge |y-x|-3 \ge \frac{9}{2}$.
Thus using  Theorem \ref{t:bhp-all}  again when $\delta_D(y) \le 3$, we get
\begin{align}\label{e:gge15}
G_D(x_1,y)
=\frac{G_D(y, x_1)}{G_D(y, x_0)}G_D(y, x_0) \le c_{13} \frac{G_D(y_1,x_1)}{G_D(y_1, x_0)}G_D(y, x_0).
\end{align}
From \eqref{e:gge14}--\eqref{e:gge15},
\begin{align}\label{e:gge16}
G_D(x,y) \le  c_{14} \frac{G_D(x_1,y_1)}{G_D(x_1, x_0)}G_D(x_0, x)  \frac{G_D(y, x_0)}{G_D(y_1, x_0)}.
\end{align}
Combining \eqref{e:gge13}
and \eqref{e:gge16} and using Step 2, we conclude that
$$
\frac{G_D(x,y)}{G_D(x,w)} \le c_{15}\frac{G_D(x_1,y_1)}{G_D(x_1,w)}
\frac{G_{D}(x_0,y)}{G_D(x_0,y_1)}
= c_{15} \left(\frac{G_D(x_1,y_1)}{G_D(x_1,w)} \frac{G_D(x_0,w)}{G_D(x_0,y_1)}  \right)
 \frac{G_{D}(x_0,y)}{G_{D}(x_0,w)} \le c_{16}\frac{G_{D}(x_0,y)}{G_{D}(x_0,w)}.
$$
\qed

For $E\subset D$ and $n\ge 1$, define
$
E^n=E\cap \{x\in D:\, 2^n\le |x|<2^{n+1}\}.
$
\begin{prop}\label{p:minthin-criterion-1-infty}
The set $E$ is minimally thin in $D$ at infinity with respect to $X$  if and only if
\begin{align}\label{e:p:m1}
\sum_{n=1}^{\infty}R^{E^n}_{M_D(\cdot, \infty)}(x_0)<\infty.\end{align}
\end{prop}
\pf
Assume
\eqref{e:p:m1} holds.
Then there
exists $n_0\in \N$ such that
$\sum_{n_0}^{\infty}R^{E^n}_{M_D(\cdot, \infty)}(x_0)<\frac12 M_D(x_0,\infty)\, .$
Let $B=\overline{B}(0,2^{n_0})$. Then $A:=B^c\cap E\subset \cup_{n_0}^{\infty} E^n$.
Therefore, $R^E_{M_D(\cdot, \infty)}(x_0)\le \sum_{n_0}^{\infty}R^{E^n}_{M_D(\cdot, \infty)}(x_0)
<\frac12 M_D(x_0,\infty)$, implying $\wh{R}^E_{M_D(\cdot, \infty)}<\frac12 M_D(x_0,\infty)$.
Hence, there exists an excessive function $u$ such that $u\ge M_D(\cdot, \infty)$ on $E$ and
$u(x_0)<\frac12 M_D(x_0,\infty)$. By the Riesz decomposition, $s=G_D\mu +M_D\nu$.
Hence, $M_D \nu(x_0)=\int_{\partial_M D} M_D(x_0,z)\, \nu(dz)< \frac12 M_D(x_0,\infty)$
implying that $\nu(\{\infty\})<\frac12$.  Therefore,
$$
\liminf_{x\to \infty, x\in E} \frac{u(x)}{M_D(x,\infty)}\ge 1 >\frac12 >\nu(\{\infty\})\, .
$$
By Proposition \ref{p:thinness} this means that $E$ is minimally thin in $D$ at infinity with respect to $X$.

Conversely, suppose that $E$ is minimally thin in $D$ at infinity with respect to $X$. By
Proposition \ref{p:thinness}, there exists a potential $u$ such that
$$
\liminf_{x\to \infty, x\in E} \frac{u(x)}{M_D(x,\infty)}=+\infty\, .
$$
Without loss of generality, we may assume that
$
u(x_0)\le (2C_7)^{-1},
$
and $C_7$ is the constant from Lemma \ref{l:Gfeinfinty}.
There exists $n_1\in \N$ with $n_1 \ge 10$
such that $u(x)>M_D(x,\infty)$ for all $x\in E\cap
\overline{B}(0,2^{n_1})^c$. Thus, $E\subset B(0,2^{n_1})\cap \{u>M_D(\cdot,\infty)\}$.
For $n\ge n_1$ define
$$
F_n=\{x\in D:\, 2^{n}\le|x|<2^{n+1}, u(x)>M_D(x,\infty)\}\qquad \textrm{and}\qquad
F=\bigcup_{n_1}^{\infty}F_n\, .
$$
Let $x\in E^n$. Since $|x|>2^{n_1}$,
we have that $u(x)>M_D(x,\infty)$ and thus $x\in F_n$. This shows that $E^n\subset F_n$,
$n\ge n_1$. Therefore, it suffices to show that
$\sum_{n_1}^{\infty}R^{F_n}_{M_D(\cdot, \infty)}(x_0)<\infty$.
Since $u>M_D(\cdot, \infty)$ on $F$, it follows that $R^{F_n}_{M_D(\cdot, \infty)}(x_0)\le u(x_0)\le c$.

Let $i\in\{1,2,3\}$. For every $n\in \N$, let
$
U_n=F_{n_1+3n+i}\, .
$
Since $i\in \{1,2,3\}$ is arbitrary, it suffices to show that $\sum_{n=1}^{\infty}
R^{U_n}_{M_D(\cdot,\infty)}(x_0)<\infty$. Let $U=\bigcup_{n=1}^{\infty}U_n$.
Then $U\subset F$ and thus $R^{U}_{M_D(\cdot, \infty)}(x_0)\le u(x_0)$.
Note that since $U$ is open, $\wh{R}^U_{M_D(\cdot, \infty)}=
R^U_{M_D(\cdot, \infty)}$ (see \cite[page 205]{BH}).  Since $u$ is a
potential, the same holds for $\wh{R}^U_{M_D(\cdot, \infty)}$, hence
there exists a measure $\mu$ such that $R^U_{M_D(\cdot, \infty)}=G_D\mu$.
Moreover, since $R^U_{M_D(\cdot, \infty)}$ is harmonic on $\overline{U}^c$
(cf.\cite[III.2.5]{BH}), $\mu(U^c)=0$. Let $\mu_n:=\mu_{|\overline{U}_n}$.
Since $\overline{U}_n$ are pairwise disjoint,
$$
\mu=\sum_{n=1}^\infty\mu_n \qquad \textrm{and} \qquad G_D\mu=
\sum_{n=1}^{\infty}G_D \mu_n\, .
$$
Fix $n\in \N$ and consider $l\in \N$, $x\in U_n$, $y\in \overline{U}_l$. If $l>n$,  then
$
|x|<2^{n_1+3n+i+1}\le 2^{1-3(l-n)}|y|\le \frac14 |y|,
$
and hence $|x-y|\ge |y|-|x|\ge \frac34 |y|$.
If $l<n$, then analogously,
$|y|\le \frac14 |x|$, hence $|x-y|\ge |x|-|y|\ge \frac34 |x|\ge \frac34 |y| $. Thus, in both cases,
$
|x-y|\ge \frac34 |y|$ for every $x\in U_n$ and  $y\in \overline{U_l}$, $l\neq n.$

Define $\mu_n'=\mu-\mu_n$ and let $x\in U_n$. By Lemma \ref{l:Gfeinfinty}
\begin{eqnarray*}
&&G_D\mu_n'(x)
= \int_D G_D(x,y)\, \mu_n'(dy)
\le C_7 M_D(x,\infty) \int_D  G_D(x_0,y)\, \mu_n'(dy) \\
&&\le C_5 M_D(x,\infty) G_D \mu(x_0)\le C_7 M_D(x,\infty) u(x_0) <\frac12 M_D(x,\infty)\, .\nonumber
\end{eqnarray*}
Since $G_D\mu_n+G_D\mu_n'=G_D\mu=R^U_{M_D(\cdot, \infty)}=M_D(\cdot, \infty)$
on $U$, it follows that $G_D\mu_n=M_D(\cdot, \infty)-G_D\mu_n'\ge M_D(\cdot, \infty)-
\frac12 M_D(\cdot, \infty)=\frac12 M_D(\cdot, \infty)$ on $U_n$. This implies that
$G_D\mu_n \ge \frac12 R^{U_n}_{M_D(\cdot, \infty)}$. Finally,
$$
\sum_{n=1}^{\infty}R^{U_n}_{M_D(\cdot, \infty)}(x_0)\le 2\sum_{n=1}^{\infty}G_D\mu_n(x_0)=2G_D\mu(x_0)<\infty\, .
$$
\qed

\begin{lemma}\label{l:MDe}
There exists $c>1$ such that
\begin{equation}\label{e:g-M-estimate-2}
c^{-1} G_D(x, x_0)|x|^d \le M_D(x,\infty)\le c G_D(x, x_0)|x|^d \quad x \in B(0, 30)^c \cap D\, .
\end{equation}
\end{lemma}
\pf
\noindent {\it Step 1}.
Assume $\delta_D(x) \ge 3$ and $|x| \ge 10$.
For $w \in B(0, 10|x|)^c \cap D$ with $\delta_D(w) \ge 3$, using the monotonicity of
Green functions  and \eqref{eGh01},
$$
\frac{G_D(x,w)}{G_D(x_0, w)G_D(x,x_0)} \asymp \frac{G_\H(x,w)}{G_{\H} (x_0, w)G_{\H}(x,x_0)}.
$$
Thus by \eqref{e:martin-kernel-estimate-infty},
\begin{align}\label{e:g-M-estimate-22}
\frac{M_D(x,\infty)}{G_D(x, x_0)}\asymp \frac{M_{\H}(x,\infty)}{G_{\H}(x, x_0)}
\asymp   \frac{\Phi(\delta_{\H}(x))^{1/2}}{ G_{\H}(x, x_0)}.
\end{align}
Now \eqref{e:g-M-estimate-2} follows from \eqref{e:g-M-estimate-22} and
\eqref{e:g-M-estimate-23} immediately.

\noindent {\it Step 2}.
Assume $\delta_D(x) \le 3$ and $|x| \ge 30$.
Let $x_1:= ( \wt x, |x|/3)$ so that $x, x_1 \in B(( \wt x, 0), |x|/2) \cap D$ and
$x_0 \notin B(( \wt x, 0), 2|x|/3) \cap D$.
In fact,
$
|x_0-( \wt x, 0)| \ge |( \wt x, 0)|-5 \ge |x|-|x_d|-5 \ge |x|-9 \ge2|x|/3.
$
Thus by Theorem \ref{t:bhp-all}, we have that
for $w \in B(0, 10|x|)^c \cap D$ with $\delta_D(w) \ge 3$,
$$
\frac{G_D(x,w)}{G_D(x,x_0)} \asymp
\frac{G_D(x_1,w)}{G_D(x_1,x_0)},
$$
which implies that
\begin{align}
\label{e:Gnew1}
\frac{M_D(x,\infty)}{G_D(x,x_0)}\asymp \frac{M_D(x_1,\infty)}{G_D(x_1,x_0)}.
\end{align}
Moreover, since
$|x|/3 \le |x_1| \le 2|x|$, by Step 1 and \eqref{e:Gnew1},
$$
\frac{M_D(x,\infty)}{G_D(x, x_0)}\asymp \frac{M_D(x_1,\infty)}{G_D(x_1, x_0)}
\asymp |x_1|^{-d}\asymp |x|^{-d}.
$$
We have proved the lemma.
\qed

Let $g(x):=G_D(x, x_0) \wedge C_8$
where
$C_8:=C_5(\sup_{s \ge 25} \Phi(s)s^{-d})>0$
so that, by \eqref{e:sbmG},
$g(x)=G_D(x, x_0)$ for every $x \in B(0, 30)^c \cap D$.
Lemma \ref{l:MDe} implies that  for $n \ge 5$,
$$
c_1^{-1}2^{nd} g(x)=c_1^{-1}2^{nd} G_D(x, x_0)\le M_D(x,\infty) \le c_1 2^{nd} G_D(x, x_0)
= c_1 2^{nd} g(x)\, ,\quad x\in E^n .
$$
This implies that for $n \ge 5$,
$
c_1^{-1}2^{nd} R^{E^n}_{g} \le R^{E^n}_{M_D(\cdot, \infty)} \le c_1 2^{nd} R^{E^n}_{g}\, .
$
In particular,
\begin{equation}\label{e:equivalence-1-3}
\sum_{n=1}^{\infty}R^{E^n}_{M_D(\cdot, \infty)}(x_0)<\infty \quad
\textrm{ if and only if}\quad \sum_{n=1}^{\infty}2^{nd} R^{E^n}_{g} (x_0)<\infty\, .
\end{equation}

Note that $\wh{R}^{E^n}_{g}$ is a potential, hence there exists a measure $\lambda_n$
(supported by $\overline{E}_n$) such that $\wh{R}^{E^n}_{g}=G_D \lambda_n$.
Also, $\wh{R}^{E^n}_{g}(\cdot)={g}(\cdot)=G_D(\cdot, x_0)$ on $\overline{E^n}$
for $n \ge 5$, hence
\begin{eqnarray*}
\wh{R}^{E^n}_{g}(x_0)
=\int_{\overline{E^n}}{g}(y)\, \lambda_n(dy)
=\int_{\overline{E}_n}\wh{R}^{E^n}_{g}(y)\, \lambda_n(dy)=\int_D \int_D G_D(x,y)\,
\lambda_n(dy)\, \lambda_n(dx)= \gamma_{g}(E^n)\, .
\end{eqnarray*}
We conclude from \eqref{e:equivalence-1-3} that
\begin{equation}\label{e:equivalence-1-4}
\sum_{n=1}^{\infty}R^{E^n}_{M_D(\cdot, \infty)}(x_0)<\infty \quad
\textrm{ if and only if}\quad \sum_{n=1}^{\infty}2^{nd} \gamma_{g}(E^n)<\infty\, .
\end{equation}
Thus we have proved the following Wiener-type criterion for minimal thinness.
\begin{corollary}\label{c:minthin-criterion-1-infty}
The set $E \subset D$ is minimally thin in $D$ at infinity with respect to $X$ if and only if
$
\sum_{n=1}^{\infty}2^{nd}$ $\gamma_{g}(E^n)<\infty.
$
\end{corollary}

\begin{remark}{\rm
Note that \cite[Part I, 11.3, Page 71]{AE}
has a similar criterion (attributed to Lelong-Ferrand) in case when $D$ is the
half-space $\mathbb{H}$ and $X$ is Brownian motion: $E$ is minimally thin at
 infinity if and only if $\sum_{n=1}^{\infty}2^{-nd} \gamma(E^n)<\infty$. Here
$\gamma(E^n)=\gamma_V(E^n)$ is the Green energy defined with respect to
the function $V(x)=x_d$ (and not ${g}(x)$) - see \cite[Part I, page 66]{AE}.
}
\end{remark}

Finally, we prove a version of Aikawa's criterion for minimal thinness.
\begin{prop}\label{p:aikawa-thinness-infty}
Suppose that $({\cal E}, {\cal F}_D)$ satisfies the Hardy inequality.
Let $\{Q_j\}_{j\ge 1}$ be a Whitney decomposition of $D$, $E\subset D$, and let $x_j$
denote the center of $Q_j$. Then $E$ is minimally thin in $D$ at infinity with respect
to $X$ if and only if
$$
\sum_{j: Q_j \subset B(0, 2^5)^c} \mathrm{dist}(0,Q_j)^{d} G_D(x_j, x_0)^2 \mathrm{Cap}_D(E\cap Q_j)<\infty\, .
$$
\end{prop}
\pf By Corollary  \ref{c:minthin-criterion-1-infty}, $E$ is minimally thin in $D$ at infinity
with respect to $X$ if and only if $\sum_{n=1}^{\infty}2^{nd} \gamma_{g}(E^n)<\infty$.
Further, let $V_n=\{x\in \R^d: 2^{n}\le |x| <2^{n+1}\}$ so that $E^n=E\cap V_n$.
If $V_n\cap Q_j \neq \emptyset$, then $\mathrm{dist}(0,Q_j)\asymp 2^{n}$. By
the quasi-additivity of $\gamma_{g}$ (Proposition \ref{p:quasi-additivity}(2)),
\begin{eqnarray*}
\sum_{n=1}^{\infty}2^{nd} \gamma_{g}(E^n)&\asymp &
\sum_{n=1}^{\infty}2^{nd} \sum_{ j: Q_j \subset B(0, 2^5)^c}
\gamma_{g}(E^n\cap Q_j)\\
&\asymp &\sum_{ j: Q_j \subset B(0, 2^5)^c}
\sum_{n, V_n\cap Q_j\neq \emptyset} \mathrm{dist}(0,Q_j)^{d}
\gamma_{g}(E^n\cap Q_j)\\
&=&\sum_{ j: Q_j \subset B(0, 2^5)^c} \mathrm{dist}(0,Q_j)^{d}
\sum_{n, V_n\cap Q_j\neq \emptyset}\gamma_{g}(E^n\cap Q_j)\\
&\asymp &\sum_{ j: Q_j \subset B(0, 2^5)^c} \mathrm{dist}(0,Q_j)^{d}
\gamma_{g}(E\cap Q_j)\, .
\end{eqnarray*}
For the last line we argue as follows: One inequality is the subadditivity of capacity.
For the other note that there exists $N\in \N$ such that for every $Q_j$,
$\sum_{n, V_n\cap Q_j\neq \emptyset}1=\sum_n 1_{V_n\cap Q_j}\le N$. Hence,
$\sum_{n, V_n\cap Q_j\neq \emptyset}\gamma_{g}(E\cap V_n\cap Q_j)\le
\sum_{n, V_n\cap Q_j\neq \emptyset}\gamma_{g}(E \cap Q_j)\le N \gamma_{g}(E \cap Q_j)$.

Finally, since $g$ satisfies the scale invariant inequality \eqref{e:u-scale-inv-har},
it follows from Lemma \ref{l:comp-gamma-cap}(2) that
$\gamma_{g}(E\cap Q_j) \asymp g(x_j)^2 \cp_D(E\cap Q_j)$. Since $g(x)=G_D(x, x_0)$
for $x\in B(0,30)^c\cap D$,
 this completes the proof.
\qed

\noindent
{\bf Proof of Theorem \ref{t:dahlberg-infinity}:} Assume that $E$ is minimally thin at $\infty$.
By Proposition \ref{p:aikawa-thinness-infty},
$$
\sum_{ j: Q_j \subset B(0, 2^5)^c} \mathrm{dist}(0,Q_j)^{d} G_D(x_j, x_0)^2
 \mathrm{Cap}_D(E\cap Q_j)<\infty\, .
$$
First note that
$$
\mathrm{Cap}_D(E\cap Q_j)\ge c_1 \sigma_1(E\cap Q_j)=c_1\int_E \1_{Q_j}(x)
\Psi(\delta_D(x)^{-1})\, dx\, .
$$
Next,
$\mathrm{dist}(0,Q_j)\asymp |x|$ for $x\in Q_j$. Therefore,
\begin{eqnarray*}
&&\sum_{ j: Q_j \subset B(0, 2^5)^c} \mathrm{dist}(0,Q_j)^{d} G_D(x_j, x_0)^2
\mathrm{Cap}_D(E\cap Q_j)\\
 &\ge &
c_2 \int_E |x|^{d} G_D(x, x_0)^2\Psi(\delta_D(x)^{-1}) \sum_{ j: Q_j
\subset B(0, 2^5)^c} \1_{Q_j}(x)\, dx\\
&  \ge  & c_2 \int_{E \cap B(0, 2^8)^c} |x|^{d} G_D(x, x_0)^2\Psi(\delta_D(x)^{-1}) \, dx\, .
\end{eqnarray*}
Conversely, assume that $E$ is the union of a subfamily of Whitney cubes. Then $E\cap Q_j$ is
either empty or equal to $Q_j$. Since $\mathrm{Cap}_D(Q_j)\asymp \sigma_1(Q_j)=
\int_{Q_j} \Psi(\delta_D(x)^{-1})\, dx$, we can reverse the first inequality in the display
above to conclude that
$$
\sum_{ j: Q_j \subset B(0, 2^5)^c}
{\rm dist}(0,Q_j)^{-d} G_D(x_j, x_0)^2 \mathrm{Cap}_D(E\cap Q_j)
\le c_3 \int_{E \cap B(0,2^5)^c} |x|^{d} G_D(x, x_0)^2\Psi(\delta_D(x)^{-1}) \, dx\, .
$$
\qed

\noindent
{\bf Proof of Corollary \ref{c:dahlberg-infinity}:}
By integrating the heat kernel estimates in \cite[Theorem 5.8]{BGR3}, one can easily get that,
for $x \in B(0, 10)^c \cap D$,
\begin{align}\label{e:newG}
G_D(x, x_0)\asymp \frac{\Phi(|x-x_0|)}{|x-x_0|^d}\left(1\wedge
\frac{\Phi(\delta_D(x))}{\Phi(|x-x_0|)}\right)^{1/2}
\left(1\wedge\frac{\Phi(\delta_D(x_0))}{\Phi(|x-x_0|)}\right)^{1/2}
\asymp
\frac{\Phi(\delta_D(x))^{1/2}}{|x|^d},
\end{align}
(see the proof of \cite[Theorem 7.3(iv)]{CKS}).
Thus the corollary immediately follows from this and Theorem \ref{t:dahlberg-infinity}.
\qed

\begin{remark}\label{r:martin-kernel-estimate}
{\rm Note that by using
\eqref{e:newG} we have
the following sharp two-sided Martin function estimates for half-space-like $C^{1, 1}$ open set $D$:
for every $z \in \partial D$,
\begin{align}\label{e:MKC11}
M_D(x, z) \asymp
\frac{\Phi(\delta_D(x))^{1/2}|x_0-z|^d}{|x-z|^d}
\end{align}
and
\begin{align}\label{e:MKC12}
M_D(x,\infty) \asymp \Phi(\delta_{D}(x))^{1/2}\, .
\end{align}
}
\end{remark}

%%%%%%%%%%%%%%%%%%%%%%%%%%%%%%%%%%%%%%%%%%%%%%%%%
%%%%%%%%%%%%         Minimal thinness of a set under the graph of a function      %%%%%%%%%%%%%%%%%%%%%%%%%%%%
%%%%%%%%%%%%%%%%%%%%%%%%%%%%%%%%%%%%%%%%%%%%%%%
\section{Minimal thinness of a set under the graph of a function}\label{graph}

In this section,
we will study minimal thinness of a set below the graph of a Lipschitz function,
both for finite and infinite boundary points. We start by recalling
Burdzy's result, cf.~\cite{Bur, Gar}.

Let $f:\R^{d-1}\to [0,\infty)$ be a Lipschitz function. The set $A=\{x=(\wt{x},x_d)
\in \H:\, 0<x_d\le f(\wt{x})\}$ is minimally thin in $\H$ with respect to Brownian
motion at $z=0$ if and only if
\begin{equation}\label{c:criterion2f}
\int_{\{|\wt{x}|<1\}}f(\wt{x})|\wt{x}|^{-d}\, d\wt{x} <\infty\, .
\end{equation}
It is shown recently in \cite{KSV6} that the same criterion for minimal thinness is true for
the subordinate Brownian motions studied there.
By using Corollary \ref{c:dahlberg} one can follow the proof of \cite[Theorem 4.4]{KSV6}
(cf.~also the proof of Theorem \ref{t:main2} below)
and show that the following Burdzy's criterion for minimal thinness holds.

\begin{prop}\label{c:dahlberg2}
Assume that
$D:=\{x=(\wt{x},x_d)\in \R^d:\, x_d> h(\wt{x})\}$ is the domain above the graph of a
bounded $C^{1,1}$ function $h$
and that $f:\R^{d-1}\to [0,\infty)$ is a Lipschitz function.
Suppose either  $h \equiv 0$ or $\gamma=1$ in \eqref{e:psi1}.
Then
the set $A:=\{x=(\wt{x},x_d)\in D:\,
 h(\wt{x})<x_d\le f(\wt{x})+ h(\wt{x})\}$ is minimally thin in $D$ at $z=
(\wt{0}, h(\wt{0})) $ with respect to
$X$
 if and only if
\eqref{c:criterion2f} holds
\end{prop}

We omit the proof and concentrate on a similar question for minimal thinness at infinity.

\begin{thm}\label{t:main2}
Suppose that  {\rm ({\bf H2})} holds and $\gamma=1$ in \eqref{e:psi1}, and
let $D=\{x=(\wt{x},x_d)\in \R^d:\, x_d> h(\wt{x})\}$ be the domain above the graph of a bounded $C^{1,1}$ function $h$.
Let $f:\R^{d-1}\to [0,\infty)$ be a Lipschitz function.
Then the set $A:=\{x=(\wt{x},x_d)\in \R^d:\, h(\wt{x}) <x_d\le f(\wt{x})+ h(\wt{x})\}$
is minimally thin in $D$ at infinity with respect to $X$  if and only if
\begin{equation}\label{c:criterion2}
\int_{\{|\wt{x}|>1\}}f(\wt{x})|\wt{x}|^{-d}\, d\wt{x} <\infty.
\end{equation}
\end{thm}

\pf Without loss of generality we may assume that $f(\wt{x})=0$ for $|\wt{x}|\le 1$,
$h(\wt 0)=0$ and $A=\{x=(\wt{x},x_d)\in \R^d:\, |\wt{x}|>1,  h(\wt{x}) <x_d\le f(\wt{x})+ h(\wt{x})\}$
 We first note that by the Lipschitz continuity of $f$ and boundedness of $h$, it follows that $|\wt{x}|\le |x|\le
 c_1|\wt{x}|$ for $x=(\wt{x},x_d)\in A$. Hence by Fubini's theorem we have
\begin{equation}\label{e:main2-1}
\int_A |x|^{-d}\, dx =\int_{|\wt{x}|>1}d\wt{x}\int \1_A(\wt{x},x_d)|x|^{-d}\, dx_d \asymp  \int_{|\wt{x}|>1}|\wt{x}|^{-d}\, d\wt{x} \int_{h(\wt{x})}^{f(\wt{x})+h(\wt{x})}dx_d=\int_{|\wt{x}|>1}f(\wt{x}) |\wt{x}|^{-d}\, d\wt{x} \, .
\end{equation}
It follows from Corollary \ref{c:dahlberg-infinity}(i) that if $A$ is minimally thin in $D$ at infinity, then \eqref{c:criterion2} holds true.

For the converse,
let $\{Q_j\}$ be a Whitney decomposition of $D$ and define $E= \cup_{Q_j\cap A\neq \emptyset}Q_j  $;
clearly $A\subset E$.
Let $Q_j^*$ be the interior of the double of $Q_j$ and
note that $\{Q_j^*\}$ has bounded multiplicity, say $N$. Moreover, if $Q_j\cap A\neq \emptyset$,
then by the Lipschitz continuity of $f$ and $h$ we have $|Q_j^*\cap A|\asymp |Q_j|$.
Moreover,
for $x \in Q_j^*$ we have
$|x|  \asymp \text{dist}(0, Q_j).$
Therefore
\begin{align}\label{e:main2-2}
&\int_A |x|^{-d}\, dx \le \int_{E }|x|^{-d}\, dx=\sum_{Q_j\cap A\neq \emptyset} \int_{Q_j  }|x|^{-d}\, dx  \le  c_2\sum_{Q_j\cap A\neq \emptyset}  |Q_j^*\cap A|  \text{dist}(0, Q_j)^{-d}
\nn \\ & \le c_3\sum_{Q_j\cap A\neq \emptyset} \int_{Q_j^*\cap A}|x|^{-d}\, dx \le c_3N \int_A |x|^{-d}\, dx\, .
\end{align}
If \eqref{c:criterion2} holds, then \eqref{e:main2-1} and \eqref{e:main2-2} imply that
$\int_E |x|^{-d} \, dx <\infty$. Hence, by Corollary \ref{c:dahlberg-infinity}(ii), $E$ is minimally thin, and thus $A$ is also minimally thin. \qed

\begin{example}\label{exam1}
{\rm
Suppose $c>0$ and $\delta \ge 0$. By Theorem \ref{t:main2},
the set $A:=\{x=(\wt{x},x_d)\in \H:\, 0<x_d\le c|\wt{x}|^{1-\delta}\}$ is minimally thin in $\H$ at infinity with respect to $X$
if and only if $\delta>0$.
}
\end{example}

\noindent
{\bf Acknowledgements:} We are grateful to the referee for
helpful comments, and in particular for suggesting a simpler
(and better) version of Theorem \ref{t:main2} and providing its proof.

%%%%%%%%%%%%%%%%%%%%%%%%%%%%%%%%%%%%%%%%%%%%%%%%%%%%%%%%  End of omitted text  %%%%%%%%%%%%%%%%%%%%%%%%%%%%%%%%%%%%%%%%%%%%%%%
\vspace{.1in}
\begin{singlespace}

%%%%%%%%%%%%%%%%%%%%%%%%%%%%%%%%%%
%%%%%%%%%%%%%%%%%%%%%%%%%%%%%%%%%%
%%%%%%%%            References
%%%%%%%%%%%%%%%%%%%%%%%%%%%%%%%%%
%%%%%%%%%%%%%%%%%%%%%%%%%%%%%%%%%%

\end{singlespace}

\end{doublespace}
\vskip 0.1truein

\parindent=0em

{\bf Panki Kim}

Department of Mathematical Sciences and Research Institute of Mathematics,

Seoul National University, Building 27, 1 Gwanak-ro, Gwanak-gu Seoul 151-747, Republic of Korea

E-mail: \texttt{pkim@snu.ac.kr}

\bigskip

{\bf Renming Song}

Department of Mathematics, University of Illinois, Urbana, IL 61801,
USA

E-mail: \texttt{rsong@math.uiuc.edu}

\bigskip

{\bf Zoran Vondra\v{c}ek}

Department of Mathematics, University of Zagreb, Zagreb, Croatia

Email: \texttt{vondra@math.hr}
\end{document}